\documentclass[1pt]{amsart}

\usepackage{amsmath,amsthm, amscd, amssymb, amsfonts}

\usepackage{hhline}

\newcommand{\Ind}{\operatorname{Ind}}
\newcommand{\ord}{\operatorname{ord}}
\newcommand{\oper}{\#}

\DeclareMathOperator*{\prode}{\boxtimes}

\newcommand\toba{{\mathfrak B }}
\newcommand\trasp{\pi}
\newcommand{\gr}{\operatorname{gr}}
\newcommand{\trid}{\triangleright}

\newcommand\compvert[2]{\genfrac{}{}{-1pt}{0}{#1}{#2}}

\newcommand\mvert[2]{\begin{tiny}\begin{matrix}#1\vspace{-4pt}\\#2\end{matrix}\end{tiny}}
\DeclareMathOperator*{\Tim}{\times}
\newcommand{\Times}[2]{\sideset{_#1}{_#2}\Tim}

\newcommand{\cx}{{\daleth}}

\newcommand{\lu}{{\nu}}
\newcommand{\ld}{{\ell}}

\newcommand{\abofe}{{A(\lu, \ld)}}

\newcommand{\bbofe}{{B(\lu, \ld)}}

\newcommand{\proy}{\Pi}

\newcommand{\Lc}{{\mathcal L}}
\newcommand{\Y}{{\mathcal Y}}
\newcommand{\R}{{\mathcal R}}
\newcommand{\W}{{\mathcal W}}
\newcommand{\Zc}{{\mathcal Z}}
\newcommand{\daga}{{\dagger}}

\newcommand{\acts}{{\rightharpoondown}}

\newcommand{\ku}{\mathbb C}
\newcommand{\K}{{\mathcal K}}
\newcommand{\Z}{{\mathbb Z}}
\newcommand{\N}{{\mathbb N}}
\newcommand{\I}{{\mathbb I}}

\newcommand{\G}{{\mathbb G}}
\newcommand{\M}{{\mathcal M}}
\newcommand{\Q}{{\mathcal Q}}
\newcommand{\F}{{\mathcal F}}
\newcommand{\C}{{\mathcal C}}
\newcommand{\D}{{\mathcal D}}
\newcommand{\Ec}{{\mathbf E}}
\newcommand{\Ee}{{\mathcal E}}
\newcommand{\Kc}{{\mathbf K}}

\newcommand{\uno}{{\bf 1}}
\newcommand{\uv}{{\bf 1_v}}
\newcommand{\uh}{{\bf 1_h}}
\newcommand{\m}{\mathcal{M}}
\newcommand{\n}{\mathcal{N}}

\newcommand{\Dc}{{\mathbf D}}
\newcommand{\B}{{\mathcal B}}
\newcommand{\T}{{\mathcal T}}
\newcommand{\Hc}{{\mathcal H}}
\newcommand{\Vc}{{\mathcal V}}
\newcommand{\Pc}{{\mathcal P}}
\newcommand{\Oc}{{\mathcal O}}
\newcommand{\ydh}{{}^H_H\mathcal{YD}}

\newcommand{\Ss}{{\mathcal S}}
\newcommand{\Vect}{\operatorname{Vec}}
\newcommand{\End}{\operatorname{End}}
\newcommand{\Aut}{\operatorname{Aut}}
\newcommand{\Int}{\operatorname{Int}}
\newcommand{\Ext}{\operatorname{Ext}}
\newcommand\tr{\operatorname{tr}}
\newcommand\Rep{\operatorname{Rep}}

\newcommand\card{\operatorname{card}}
\newcommand{\unosigma}{{\bf 1}^{\sigma}}
\newcommand\tot{\operatorname{tot}}
\newcommand\sgn{\operatorname{sgn}}
\newcommand\ad{\operatorname{ad}}
\newcommand\Hom{\operatorname{Hom}}
\newcommand\opext{\operatorname{Opext}}
\newcommand\Tot{\operatorname{Tot}}
\newcommand\Map{\operatorname{Map}}
\newcommand{\fde}{{\triangleright}}
\newcommand{\fiz}{{\triangleleft}}
\newcommand{\wC}{\widehat{C}}
\newcommand{\la}{\langle}
\newcommand{\ra}{\rangle}
\newcommand{\Comod}{\mbox{\rm Comod\,}}
\newcommand{\Mod}{\mbox{\rm Mod\,}}
\newcommand{\Sg}{{\mathfrak S}}
\newcommand{\funcion}{\mathfrak p}
\newcommand{\tauo}{{\widehat\tau}}
\newcommand{\prin}{t}
\newcommand{\fin}{b}
\newcommand{\pri}{r}
\newcommand{\fine}{l}
\newcommand\rh{\sim_{H}}
\newcommand\rv{\sim_{V}}
\newcommand\rd{\sim_{D}}

\newcommand\V{\operatorname{Vec}}

\theoremstyle{plain}

\renewcommand{\baselinestretch}{1.2}

\title{The character tables of centralizers in  Weyl Group of $E_8$ III}
\author{ \small Shouchuan Zhang,     Peng Wang,  Jing Cheng, Hui Yang
}
\address{ Mathematics Department of Hunan University,\newline \indent Changsha China,
410082, E-mail:\\ z9491@yahoo.com.cn }
\date{}

\begin{document}
\newtheorem{Proposition}{\quad Proposition}[section]
\newtheorem{Theorem}{\quad Theorem}
\newtheorem{Definition}[Proposition]{\quad Definition}
\newtheorem{Corollary}[Proposition]{\quad Corollary}
\newtheorem{Lemma}[Proposition]{\quad Lemma}
\newtheorem{Example}[Proposition]{\quad Example}

\maketitle \addtocounter{section}{-1}

\numberwithin{equation}{section}

\date{}


\begin {abstract} To classify the finite dimensional pointed Hopf
algebras with Weyl group $G$ of $E_8$, we obtain  the
representatives of conjugacy classes of $G$ and  all character
tables of centralizers of these representatives by means of software
GAP. In this paper we only list character table 47--64.


keywords: GAP, Hopf algebra, Weyl group, character.
\end {abstract}

\section{Introduction}\label {s0}

This article is to contribute to the classification of
finite-dimensional complex pointed Hopf algebras  with Weyl groups
of $E_8$.
 Many papers are about the classification of finite dimensional
pointed Hopf algebras, for example,  \cite{AS98b, AS02, AS00, AS05,
He06, AHS08, AG03, AFZ08, AZ07,Gr00,  Fa07, AF06, AF07, ZZC04, ZCZ08}.

In these research  ones need  the centralizers and character tables
of groups. In this paper we obtain   the representatives of
conjugacy classes of Weyl groups of $E_8$ and all
character tables of centralizers of these representatives by means
of software GAP. In this paper we only list character table  47--64.

By the Cartan-Killing classification of simple Lie algebras over complex field the Weyl groups
to
be considered are $W(A_l), (l\geq 1); $
 $W(B_l), (l \geq 2);$
 $W(C_l), (l \geq 2); $
 $W(D_l), (l\geq 4); $
 $W(E_l), (8 \geq l \geq 6); $
 $W(F_l), ( l=4 );$
$W(G_l), (l=2)$.

It is otherwise desirable to do this in view of the importance of Weyl groups in the theories of
Lie groups, Lie algebras and algebraic groups. For example, the irreducible representations of
Weyl groups were obtained by Frobenius, Schur and Young. The conjugace classes of  $W(F_4)$
were obtained by Wall
 \cite {Wa63} and its character tables were obtained by Kondo \cite {Ko65}.
The conjugace classes  and  character tables of $W(E_6),$ $W(E_7)$ and $W(E_8)$ were obtained
by  Frame \cite {Fr51}.
Carter gave a unified description of the conjugace classes of
 Weyl groups of simple Lie algebras  \cite  {Ca72}.

\section {Program}

 By using the following
program in GAP,  we  obtain  the representatives of conjugacy
classes of Weyl groups of $E_6$  and all character tables of
centralizers of these representatives.

gap$>$ L:=SimpleLieAlgebra("E",6,Rationals);;

gap$>$ R:=RootSystem(L);;

 gap$>$ W:=WeylGroup(R);Display(Order(W));

gap $>$ ccl:=ConjugacyClasses(W);;

gap$>$ q:=NrConjugacyClasses(W);; Display (q);

gap$>$ for i in [1..q] do

$>$ r:=Order(Representative(ccl[i]));Display(r);;

$>$ od; gap

$>$ s1:=Representative(ccl[1]);;cen1:=Centralizer(W,s1);;

gap$>$ cl1:=ConjugacyClasses(cen1);

gap$>$ s1:=Representative(ccl[2]);;cen1:=Centralizer(W,s1);;

 gap$>$
cl2:=ConjugacyClasses(cen1);

gap$>$ s1:=Representative(ccl[3]);;cen1:=Centralizer(W,s1);;

gap$>$ cl3:=ConjugacyClasses(cen1);

 gap$>$
s1:=Representative(ccl[4]);;cen1:=Centralizer(W,s1);;

 gap$>$
cl4:=ConjugacyClasses(cen1);

 gap$>$
s1:=Representative(ccl[5]);;cen1:=Centralizer(W,s1);;

gap$>$ cl5:=ConjugacyClasses(cen1);

 gap$>$
s1:=Representative(ccl[6]);;cen1:=Centralizer(W,s1);;

gap$>$ cl6:=ConjugacyClasses(cen1);

gap$>$ s1:=Representative(ccl[7]);;cen1:=Centralizer(W,s1);;

gap$>$ cl7:=ConjugacyClasses(cen1);

gap$>$ s1:=Representative(ccl[8]);;cen1:=Centralizer(W,s1);;

gap$>$ cl8:=ConjugacyClasses(cen1);

gap$>$ s1:=Representative(ccl[9]);;cen1:=Centralizer(W,s1);;

gap$>$ cl9:=ConjugacyClasses(cen1);

gap$>$ s1:=Representative(ccl[10]);;cen1:=Centralizer(W,s1);;

 gap$>$
cl10:=ConjugacyClasses(cen1);

gap$>$ s1:=Representative(ccl[11]);;cen1:=Centralizer(W,s1);;

gap$>$ cl11:=ConjugacyClasses(cen1);

gap$>$ s1:=Representative(ccl[12]);;cen1:=Centralizer(W,s1);;

gap$>$ cl2:=ConjugacyClasses(cen1);

gap$>$ s1:=Representative(ccl[13]);;cen1:=Centralizer(W,s1);;

 gap$>$
cl13:=ConjugacyClasses(cen1);

 gap$>$
s1:=Representative(ccl[14]);;cen1:=Centralizer(W,s1);;

gap$>$ cl14:=ConjugacyClasses(cen1);

gap$>$ s1:=Representative(ccl[15]);;cen1:=Centralizer(W,s1);;

gap$>$ cl15:=ConjugacyClasses(cen1);

gap$>$ s1:=Representative(ccl[16]);;cen1:=Centralizer(W,s1);;

gap$>$ cl16:=ConjugacyClasses(cen1);

gap$>$ s1:=Representative(ccl[17]);;cen1:=Centralizer(W,s1);;

gap$>$ cl17:=ConjugacyClasses(cen1);

gap$>$ s1:=Representative(ccl[18]);;cen1:=Centralizer(W,s1);;

gap$>$ cl18:=ConjugacyClasses(cen1);

gap$>$ s1:=Representative(ccl[19]);;cen1:=Centralizer(W,s1);;

gap$>$ cl19:=ConjugacyClasses(cen1);

gap$>$ s1:=Representative(ccl[20]);;cen1:=Centralizer(W,s1);;

gap$>$ cl20:=ConjugacyClasses(cen1);

gap$>$ s1:=Representative(ccl[21]);;cen1:=Centralizer(W,s1);;

 gap$>$
cl21:=ConjugacyClasses(cen1);

gap$>$ s1:=Representative(ccl[22]);;cen1:=Centralizer(W,s1);;

 gap$>$
cl22:=ConjugacyClasses(cen1);

gap$>$ s1:=Representative(ccl[23]);;cen1:=Centralizer(W,s1);;

gap$>$ cl23:=ConjugacyClasses(cen1);

gap$>$ s1:=Representative(ccl[24]);;cen1:=Centralizer(W,s1);;

$>$ cl24:=ConjugacyClasses(cen1);

gap$>$ s1:=Representative(ccl[25]);;cen1:=Centralizer(W,s1);

gap$>$ cl25:=ConjugacyClasses(cen1);

gap$>$ for i in [1..q] do

$>$ s:=Representative(ccl[i]);;cen:=Centralizer(W,s);;

$>$ char:=CharacterTable(cen);;Display (cen);Display(char);

 $>$ od;

The programs for Weyl groups of $E_7$, $E_8$, $F_4$ and $ G_2$ are
similar.

{\small

\section {$E_8$}

The generators of $G^{s_{47}}$ are:\\
$\left(
\begin{array}{cccccccc} -2&-1&1&0& 1&-1&1&0 \\
-2&-2&2&0&1&-2&2&0 \\ -3&-3&2&0&2&-2&
2&0 \\
        -4&-4&4&-1&3&-3&3&0 \\ -3&-3&3&-1&2&-2&3&0 \\
        -2&-3&2&0&1&-1&2&0 \\ -2&-2&2&0&0&0&1&0 \\ -1&-1&1&0&0&0&0&1 \end{array}\right)$,
  $\left(\begin{array}{cccccccc} 1&0&0&0&0&0&0&0 \\ 0&1&0&0&0&0&0&0 \\ 0&0&1&0&0&0&0&0 \\
        0&0&0&1&0&0&0&0 \\ 0&0&0&0&1&0&0&0 \\ 0&0&0&0&0&1&0&0 \\
        0&0&0&0&0&0&1&0 \\ 0&0&0&0&0&0&1&-1 \end{array}\right)$,
  $\left(\begin{array}{cccccccc} 1&0&0&0&0&0&0&0 \\ 1&0&-1&1&0&-1&1&0 \\ 1&-1&0&1&0&-1&1&0 \\
        1&-1&-1&2&0&-1&1&0 \\ 1&-1&-1&1&1&-1&1&0 \\ 1&-1&-1&1&0&0&1&0 \\
        0&0&0&0&0&0&1&0 \\ 0&0&0&0&0&0&0&1 \end{array}\right)$,
  $\left(\begin{array}{cccccccc} 1&0&-1&0&0&2&-2&0 \\ 2&0&-3&1&0&2&-2&0 \\ 2&0&-3&1&0&3&-3&0 \\
        4&-1&-5&2&0&4&-4&0 \\ 3&-1&-4&1&1&3&-3&0 \\ 2&-1&-3&1&0&3&-2&0 \\
        2&0&-2&0&0&2&-1&0 \\ 1&0&-1&0&0&1&-1&1 \end{array}\right)$,
  $\left(\begin{array}{cccccccc} 1&1&0&-1&1&-1&0&1 \\ 0&2&0&-1&1&-1&0&1 \\ 0&2&1&-2&2&-2&0&2 \\
        0&3&0&-2&3&-3&0&3 \\ 0&3&0&-2&2&-2&0&3 \\ 0&2&0&-1&1&-1&-1&3 \\
        0&1&0&0&0&-1&0&2 \\ 0&0&0&0&0&0&0&1 \end{array}\right)$,
  $\left(\begin{array}{cccccccc} -1&0&0&2&-2&0&0&0 \\ -1&0&1&2&-3&1&-1&0 \\ -2&0&1&3&-3&0&-1&0 \\
        -2&0&2&4&-5&0&-1&0 \\ -2&0&2&3&-4&0&-1&0 \\ -1&1&1&2&-3&0&-1&0 \\
        0&0&0&2&-2&0&-1&0 \\ 0&0&0&1&-1&0&-1&1 \end{array}\right)$.

The representatives of conjugacy classes of   $G^{s_{47}}$ are:\\
$\left(\begin{array}{cccccccc} 1&0&0&0&0&0&0&0 \\ 0&1&0&0&
0&0&0&0 \\ 0&0&1&0&0&0&0&0 \\
        0&0&0&1&0&0&0&0 \\ 0&0&0&0&1&0&0&0 \\ 0&0&0&0&0&1&0&0 \\
        0&0&0&0&0&0&1&0 \\ 0&0&0&0&0&0&0&1 \end{array}\right)$,
  $\left(\begin{array}{cccccccc} 0&-1&1&-1&0&1&0&0 \\ -1&-1&1&-1&0&2&-1&0 \\ 0&-1&1&-2&0&3&-1&0 \\
        -1&-2&1&-2&0&4&-1&0 \\ -1&-2&1&-1&-1&4&-1&0 \\
        -1&-2&1&-1&0&3&-1&0 \\ 0&-2&0&0&0&2&-1&0 \\ 0&-1&0&0&0&1&-1&1 \end{array}\right)$,
 $\left(\begin{array}{cccccccc} -2&-1&1&0&1&-1&1&0 \\ -2&-2&2&0&1&-2&2&0 \\
        -3&-3&2&0&2&-2&2&0 \\ -4&-4&4&-1&3&-3&3&0 \\
        -3&-3&3&-1&2&-2&3&0 \\ -2&-3&2&0&1&-1&2&0 \\ -2&-2&2&0&0&0&1&0 \\
        -1&-1&1&0&0&0&0&1 \end{array}\right)$,
  $\left(\begin{array}{cccccccc} 1&0&-2&1&-1&2&-1&0 \\ 3&0&-3&1&-1&2&-1&0 \\ 3&0&-4&2&-2&3&-1&0 \\
        5&0&-5&2&-3&5&-2&0 \\ 4&-1&-4&2&-2&4&-2&0 \\
        3&-1&-3&1&-1&3&-1&0 \\ 2&0&-2&0&0&2&-1&0 \\ 1&0&-1&0&0&1&-1&1 \end{array}\right)$,
  $\left(\begin{array}{cccccccc} 1&0&0&0&0&0&0&0 \\ 0&1&0&0&0&0&0&0 \\ 0&0&1&0&0&0&0&0 \\
        0&0&0&1&0&0&0&0 \\ 0&0&0&0&1&0&0&0 \\ 0&0&0&0&0&1&0&0 \\
        0&0&0&0&0&0&1&0 \\ 0&0&0&0&0&0&1&-1 \end{array}\right)$,
  $\left(\begin{array}{cccccccc} 0&-1&1&-1&0&1&0&0 \\ -1&-1&1&-1&0&2&-1&0 \\ 0&-1&1&-2&0&3&-1&0 \\
        -1&-2&1&-2&0&4&-1&0 \\ -1&-2&1&-1&-1&4&-1&0 \\
        -1&-2&1&-1&0&3&-1&0 \\ 0&-2&0&0&0&2&-1&0 \\ 0&-1&0&0&0&1&0&-1 \end{array}\right)$,
 $\left(\begin{array}{cccccccc} -2&-1&1&0&1&-1&1&0 \\ -2&-2&2&0&1&-2&2&0 \\
        -3&-3&2&0&2&-2&2&0 \\ -4&-4&4&-1&3&-3&3&0 \\
        -3&-3&3&-1&2&-2&3&0 \\ -2&-3&2&0&1&-1&2&0 \\ -2&-2&2&0&0&0&1&0 \\
        -1&-1&1&0&0&0&1&-1 \end{array}\right)$,
  $\left(\begin{array}{cccccccc} 1&0&-2&1&-1&2&-1&0 \\ 3&0&-3&1&-1&2&-1&0 \\ 3&0&-4&2&-2&3&-1&0 \\
        5&0&-5&2&-3&5&-2&0 \\ 4&-1&-4&2&-2&4&-2&0 \\
        3&-1&-3&1&-1&3&-1&0 \\ 2&0&-2&0&0&2&-1&0 \\ 1&0&-1&0&0&1&0&-1 \end{array}\right)$,
  $\left(\begin{array}{cccccccc} 1&1&0&0&0&-1&0&0 \\ 0&2&0&0&0&-1&0&0 \\ 0&2&1&0&0&-2&0&0 \\
        0&3&0&1&0&-3&0&0 \\ 0&3&0&0&1&-3&0&0 \\ 0&3&0&0&0&-2&0&0 \\
        0&2&0&0&0&-2&1&0 \\ 0&1&0&0&0&-1&1&-1 \end{array}\right)$,
  $\left(\begin{array}{cccccccc} 0&0&1&-1&0&0&0&0 \\ -1&0&1&-1&0&1&-1&0 \\ 0&1&1&-2&0&1&-1&0 \\
        -1&1&1&-2&0&1&-1&0 \\ -1&1&1&-1&-1&1&-1&0 \\
        -1&1&1&-1&0&0&-1&0 \\ 0&0&0&0&0&0&-1&0 \\ 0&0&0&0&0&0&0&-1 \end{array}\right)$,
  $\left(\begin{array}{cccccccc} -2&0&1&0&1&-2&1&0 \\ -2&-1&2&0&1&-3&2&0 \\ -3&-1&2&0&2&-4&2&0 \\
        -4&-1&4&-1&3&-6&3&0 \\ -3&0&3&-1&2&-5&3&0 \\ -2&0&2&0&1&-4&2&0 \\
        -2&0&2&0&0&-2&1&0 \\ -1&0&1&0&0&-1&1&-1 \end{array}\right)$,
  $\left(\begin{array}{cccccccc} 1&1&-2&1&-1&1&-1&0 \\ 3&1&-3&1&-1&1&-1&0 \\ 3&2&-4&2&-2&1&-1&0 \\
        5&3&-5&2&-3&2&-2&0 \\ 4&2&-4&2&-2&1&-2&0 \\ 3&2&-3&1&-1&0&-1&0 \\
        2&2&-2&0&0&0&-1&0 \\ 1&1&-1&0&0&0&0&-1 \end{array}\right)$,
  $\left(\begin{array}{cccccccc} 1&0&0&0&0&0&0&0 \\ 1&0&-1&1&0&-1&1&0 \\ 1&-1&0&1&0&-1&1&0 \\
        1&-1&-1&2&0&-1&1&0 \\ 1&-1&-1&1&1&-1&1&0 \\ 1&-1&-1&1&0&0&1&0 \\
        0&0&0&0&0&0&1&0 \\ 0&0&0&0&0&0&0&1 \end{array}\right)$,
  $\left(\begin{array}{cccccccc} 0&-1&1&-1&0&1&0&0 \\ 0&-2&0&0&0&1&0&0 \\ 1&-2&0&-1&0&2&0&0 \\
        0&-3&0&-1&0&3&0&0 \\ 0&-3&0&0&-1&3&0&0 \\ 0&-3&0&0&0&2&0&0 \\
        0&-2&0&0&0&2&-1&0 \\ 0&-1&0&0&0&1&-1&1 \end{array}\right)$,
  $\left(\begin{array}{cccccccc} -2&-1&1&0&1&-1&1&0 \\ -3&-1&3&-1&1&-1&1&0 \\ -4&-2&3&-1&2&-1&1&0 \\
        -5&-3&5&-2&3&-2&2&0 \\ -4&-2&4&-2&2&-1&2&0 \\
        -3&-2&3&-1&1&0&1&0 \\ -2&-2&2&0&0&0&1&0 \\ -1&-1&1&0&0&0&0&1 \end{array}\right)$,
 $\left(\begin{array}{cccccccc} 1&0&-2&1&-1&2&-1&0 \\ 2&1&-2&0&-1&3&-2&0 \\
        2&1&-3&1&-2&4&-2&0 \\ 4&1&-4&1&-3&6&-3&0 \\ 3&0&-3&1&-2&5&-3&0 \\
        2&0&-2&0&-1&4&-2&0 \\ 2&0&-2&0&0&2&-1&0 \\ 1&0&-1&0&0&1&-1&1 \end{array}\right)$,
  $\left(\begin{array}{cccccccc} 1&0&0&0&0&0&0&0 \\ 1&0&-1&1&0&-1&1&0 \\ 1&-1&0&1&0&-1&1&0 \\
        1&-1&-1&2&0&-1&1&0 \\ 1&-1&-1&1&1&-1&1&0 \\ 1&-1&-1&1&0&0&1&0 \\
        0&0&0&0&0&0&1&0 \\ 0&0&0&0&0&0&1&-1 \end{array}\right)$,
  $\left(\begin{array}{cccccccc} 0&-1&1&-1&0&1&0&0 \\ 0&-2&0&0&0&1&0&0 \\ 1&-2&0&-1&0&2&0&0 \\
        0&-3&0&-1&0&3&0&0 \\ 0&-3&0&0&-1&3&0&0 \\ 0&-3&0&0&0&2&0&0 \\
        0&-2&0&0&0&2&-1&0 \\ 0&-1&0&0&0&1&0&-1 \end{array}\right)$,
  $\left(\begin{array}{cccccccc} -2&-1&1&0&1&-1&1&0 \\ -3&-1&3&-1&1&-1&1&0 \\ -4&-2&3&-1&2&-1&1&0 \\
        -5&-3&5&-2&3&-2&2&0 \\ -4&-2&4&-2&2&-1&2&0 \\
        -3&-2&3&-1&1&0&1&0 \\ -2&-2&2&0&0&0&1&0 \\ -1&-1&1&0&0&0&1&-1 \end{array}\right)$,
   $\left(\begin{array}{cccccccc} 1&0&-2&1&-1&2&-1&0 \\ 2&1&-2&0&-1&3&-2&0 \\
        2&1&-3&1&-2&4&-2&0 \\ 4&1&-4&1&-3&6&-3&0 \\ 3&0&-3&1&-2&5&-3&0 \\
        2&0&-2&0&-1&4&-2&0 \\ 2&0&-2&0&0&2&-1&0 \\ 1&0&-1&0&0&1&0&-1 \end{array}\right)$,
  $\left(\begin{array}{cccccccc} 1&1&0&0&0&-1&0&0 \\ 1&1&-1&1&0&-2&1&0 \\ 1&1&0&1&0&-3&1&0 \\
        1&2&-1&2&0&-4&1&0 \\ 1&2&-1&1&1&-4&1&0 \\ 1&2&-1&1&0&-3&1&0 \\
        0&2&0&0&0&-2&1&0 \\ 0&1&0&0&0&-1&1&-1 \end{array}\right)$,
  $\left(\begin{array}{cccccccc} 0&0&1&-1&0&0&0&0 \\ 0&-1&0&0&0&0&0&0 \\ 1&0&0&-1&0&0&0&0 \\
        0&0&0&-1&0&0&0&0 \\ 0&0&0&0&-1&0&0&0 \\ 0&0&0&0&0&-1&0&0 \\
        0&0&0&0&0&0&-1&0 \\ 0&0&0&0&0&0&0&-1 \end{array}\right)$,
  $\left(\begin{array}{cccccccc} -2&0&1&0&1&-2&1&0 \\ -3&0&3&-1&1&-2&1&0 \\ -4&0&3&-1&2&-3&1&0 \\
        -5&0&5&-2&3&-5&2&0 \\ -4&1&4&-2&2&-4&2&0 \\ -3&1&3&-1&1&-3&1&0 \\
        -2&0&2&0&0&-2&1&0 \\ -1&0&1&0&0&-1&1&-1 \end{array}\right)$,
  $\left(\begin{array}{cccccccc} 1&1&-2&1&-1&1&-1&0 \\ 2&2&-2&0&-1&2&-2&0 \\ 2&3&-3&1&-2&2&-2&0 \\
        4&4&-4&1&-3&3&-3&0 \\ 3&3&-3&1&-2&2&-3&0 \\ 2&3&-2&0&-1&1&-2&0 \\
        2&2&-2&0&0&0&-1&0 \\ 1&1&-1&0&0&0&0&-1 \end{array}\right)$,
  $\left(\begin{array}{cccccccc} 0&0&-1&1&0&0&0&0 \\ 1&0&-1&1&0&-1&1&0 \\ 0&-1&-1&2&0&-1&1&0 \\
        1&-1&-1&2&0&-1&1&0 \\ 1&-1&-1&1&1&-1&1&0 \\ 1&-1&-1&1&0&0&1&0 \\
        0&0&0&0&0&0&1&0 \\ 0&0&0&0&0&0&0&1 \end{array}\right)$,
  $\left(\begin{array}{cccccccc} -1&-1&0&0&0&1&0&0 \\ 0&-2&0&0&0&1&0&0 \\ 0&-2&-1&0&0&2&0&0 \\
        0&-3&0&-1&0&3&0&0 \\ 0&-3&0&0&-1&3&0&0 \\ 0&-3&0&0&0&2&0&0 \\
        0&-2&0&0&0&2&-1&0 \\ 0&-1&0&0&0&1&-1&1 \end{array}\right)$,
  $\left(\begin{array}{cccccccc} -1&-1&2&-1&1&-1&1&0 \\ -3&-1&3&-1&1&-1&1&0 \\
        -3&-2&4&-2&2&-1&1&0 \\ -5&-3&5&-2&3&-2&2&0 \\
        -4&-2&4&-2&2&-1&2&0 \\ -3&-2&3&-1&1&0&1&0 \\ -2&-2&2&0&0&0&1&0 \\
        -1&-1&1&0&0&0&0&1 \end{array}\right)$,
  $\left(\begin{array}{cccccccc} 2&0&-1&0&-1&2&-1&0 \\ 2&1&-2&0&-1&3&-2&0 \\ 3&1&-2&0&-2&4&-2&0 \\
        4&1&-4&1&-3&6&-3&0 \\ 3&0&-3&1&-2&5&-3&0 \\ 2&0&-2&0&-1&4&-2&0 \\
        2&0&-2&0&0&2&-1&0 \\ 1&0&-1&0&0&1&-1&1 \end{array}\right)$,
  $\left(\begin{array}{cccccccc} 0&0&-1&1&0&0&0&0 \\ 1&0&-1&1&0&-1&1&0 \\ 0&-1&-1&2&0&-1&1&0 \\
        1&-1&-1&2&0&-1&1&0 \\ 1&-1&-1&1&1&-1&1&0 \\ 1&-1&-1&1&0&0&1&0 \\
        0&0&0&0&0&0&1&0 \\ 0&0&0&0&0&0&1&-1 \end{array}\right)$,
  $\left(\begin{array}{cccccccc} -1&-1&0&0&0&1&0&0 \\ 0&-2&0&0&0&1&0&0 \\ 0&-2&-1&0&0&2&0&0 \\
        0&-3&0&-1&0&3&0&0 \\ 0&-3&0&0&-1&3&0&0 \\ 0&-3&0&0&0&2&0&0 \\
        0&-2&0&0&0&2&-1&0 \\ 0&-1&0&0&0&1&0&-1 \end{array}\right)$,
  $\left(\begin{array}{cccccccc} -1&-1&2&-1&1&-1&1&0 \\ -3&-1&3&-1&1&-1&1&0 \\
        -3&-2&4&-2&2&-1&1&0 \\ -5&-3&5&-2&3&-2&2&0 \\
        -4&-2&4&-2&2&-1&2&0 \\ -3&-2&3&-1&1&0&1&0 \\ -2&-2&2&0&0&0&1&0 \\
        -1&-1&1&0&0&0&1&-1 \end{array}\right)$,
  $\left(\begin{array}{cccccccc} 2&0&-1&0&-1&2&-1&0 \\ 2&1&-2&0&-1&3&-2&0 \\ 3&1&-2&0&-2&4&-2&0 \\
        4&1&-4&1&-3&6&-3&0 \\ 3&0&-3&1&-2&5&-3&0 \\ 2&0&-2&0&-1&4&-2&0 \\
        2&0&-2&0&0&2&-1&0 \\ 1&0&-1&0&0&1&0&-1 \end{array}\right)$,
  $\left(\begin{array}{cccccccc} 0&1&-1&1&0&-1&0&0 \\ 1&1&-1&1&0&-2&1&0 \\ 0&1&-1&2&0&-3&1&0 \\
        1&2&-1&2&0&-4&1&0 \\ 1&2&-1&1&1&-4&1&0 \\ 1&2&-1&1&0&-3&1&0 \\
        0&2&0&0&0&-2&1&0 \\ 0&1&0&0&0&-1&1&-1 \end{array}\right)$,
  $\left(\begin{array}{cccccccc} -1&0&0&0&0&0&0&0 \\ 0&-1&0&0&0&0&0&0 \\ 0&0&-1&0&0&0&0&0 \\
        0&0&0&-1&0&0&0&0 \\ 0&0&0&0&-1&0&0&0 \\ 0&0&0&0&0&-1&0&0 \\
        0&0&0&0&0&0&-1&0 \\ 0&0&0&0&0&0&0&-1 \end{array}\right)$,
  $\left(\begin{array}{cccccccc} -1&0&2&-1&1&-2&1&0 \\ -3&0&3&-1&1&-2&1&0 \\ -3&0&4&-2&2&-3&1&0 \\
        -5&0&5&-2&3&-5&2&0 \\ -4&1&4&-2&2&-4&2&0 \\ -3&1&3&-1&1&-3&1&0 \\
        -2&0&2&0&0&-2&1&0 \\ -1&0&1&0&0&-1&1&-1 \end{array}\right)$,
  $\left(\begin{array}{cccccccc} 2&1&-1&0&-1&1&-1&0 \\ 2&2&-2&0&-1&2&-2&0 \\ 3&3&-2&0&-2&2&-2&0 \\
        4&4&-4&1&-3&3&-3&0 \\ 3&3&-3&1&-2&2&-3&0 \\ 2&3&-2&0&-1&1&-2&0 \\
        2&2&-2&0&0&0&-1&0 \\ 1&1&-1&0&0&0&0&-1 \end{array}\right)$,
  $\left(\begin{array}{cccccccc} -1&-2&0&1&0&1&-1&0 \\ 0&-3&0&1&0&1&-1&0 \\ 0&-4&-1&2&0&2&-2&0 \\
        0&-6&0&2&0&3&-3&0 \\ 0&-5&0&1&1&2&-2&0 \\ 0&-4&0&1&0&2&-1&0 \\
        0&-2&0&0&0&2&-1&0 \\ 0&-1&0&0&0&1&-1&1 \end{array}\right)$,
  $\left(\begin{array}{cccccccc} 2&0&-1&-1&1&1&-1&0 \\ 2&1&-2&-1&1&2&-2&0 \\ 3&1&-2&-2&2&2&-2&0 \\
        4&1&-4&-2&3&3&-3&0 \\ 3&1&-3&-2&2&3&-2&0 \\ 2&0&-2&-1&1&3&-2&0 \\
        2&0&-2&0&0&2&-1&0 \\ 1&0&-1&0&0&1&-1&1 \end{array}\right)$,
  $\left(\begin{array}{cccccccc} -1&-2&0&1&0&1&-1&0 \\ 0&-3&0&1&0&1&-1&0 \\ 0&-4&-1&2&0&2&-2&0 \\
        0&-6&0&2&0&3&-3&0 \\ 0&-5&0&1&1&2&-2&0 \\ 0&-4&0&1&0&2&-1&0 \\
        0&-2&0&0&0&2&-1&0 \\ 0&-1&0&0&0&1&0&-1 \end{array}\right)$,
  $\left(\begin{array}{cccccccc} 2&0&-1&-1&1&1&-1&0 \\ 2&1&-2&-1&1&2&-2&0 \\ 3&1&-2&-2&2&2&-2&0 \\
        4&1&-4&-2&3&3&-3&0 \\ 3&1&-3&-2&2&3&-2&0 \\ 2&0&-2&-1&1&3&-2&0 \\
        2&0&-2&0&0&2&-1&0 \\ 1&0&-1&0&0&1&0&-1 \end{array}\right)$,
  $\left(\begin{array}{cccccccc} -1&-1&0&1&0&0&-1&0 \\ 0&-2&0&1&0&0&-1&0 \\ 0&-2&-1&2&0&0&-2&0 \\
        0&-3&0&2&0&0&-3&0 \\ 0&-2&0&1&1&-1&-2&0 \\ 0&-1&0&1&0&-1&-1&0 \\
        0&0&0&0&0&0&-1&0 \\ 0&0&0&0&0&0&0&-1 \end{array}\right)$,
  $\left(\begin{array}{cccccccc} 2&1&-1&-1&1&0&-1&0 \\ 2&2&-2&-1&1&1&-2&0 \\ 3&3&-2&-2&2&0&-2&0 \\
        4&4&-4&-2&3&0&-3&0 \\ 3&4&-3&-2&2&0&-2&0 \\ 2&3&-2&-1&1&0&-2&0 \\
        2&2&-2&0&0&0&-1&0 \\ 1&1&-1&0&0&0&0&-1 \end{array}\right)$,
  $\left(\begin{array}{cccccccc} -1&-2&0&1&0&1&-1&0 \\ -1&-2&1&0&0&2&-2&0 \\ -1&-3&0&1&0&3&-3&0 \\
        -1&-5&1&1&0&4&-4&0 \\ -1&-4&1&0&1&3&-3&0 \\ -1&-3&1&0&0&3&-2&0 \\
        0&-2&0&0&0&2&-1&0 \\ 0&-1&0&0&0&1&-1&1 \end{array}\right)$,
  $\left(\begin{array}{cccccccc} 2&0&-1&-1&1&1&-1&0 \\ 3&0&-3&0&1&1&-1&0 \\ 4&0&-3&-1&2&1&-1&0 \\
        5&0&-5&-1&3&2&-2&0 \\ 4&0&-4&-1&2&2&-1&0 \\ 3&-1&-3&0&1&2&-1&0 \\
        2&0&-2&0&0&2&-1&0 \\ 1&0&-1&0&0&1&-1&1 \end{array}\right)$,
  $\left(\begin{array}{cccccccc} -1&-2&0&1&0&1&-1&0 \\ -1&-2&1&0&0&2&-2&0 \\ -1&-3&0&1&0&3&-3&0 \\
        -1&-5&1&1&0&4&-4&0 \\ -1&-4&1&0&1&3&-3&0 \\ -1&-3&1&0&0&3&-2&0 \\
        0&-2&0&0&0&2&-1&0 \\ 0&-1&0&0&0&1&0&-1 \end{array}\right)$,
  $\left(\begin{array}{cccccccc} 2&0&-1&-1&1&1&-1&0 \\ 3&0&-3&0&1&1&-1&0 \\ 4&0&-3&-1&2&1&-1&0 \\
        5&0&-5&-1&3&2&-2&0 \\ 4&0&-4&-1&2&2&-1&0 \\ 3&-1&-3&0&1&2&-1&0 \\
        2&0&-2&0&0&2&-1&0 \\ 1&0&-1&0&0&1&0&-1 \end{array}\right)$,
  $\left(\begin{array}{cccccccc} -1&-1&0&1&0&0&-1&0 \\ -1&-1&1&0&0&1&-2&0 \\ -1&-1&0&1&0&1&-3&0 \\
        -1&-2&1&1&0&1&-4&0 \\ -1&-1&1&0&1&0&-3&0 \\ -1&0&1&0&0&0&-2&0 \\
        0&0&0&0&0&0&-1&0 \\ 0&0&0&0&0&0&0&-1 \end{array}\right)$,
  $\left(\begin{array}{cccccccc} 2&1&-1&-1&1&0&-1&0 \\ 3&1&-3&0&1&0&-1&0 \\ 4&2&-3&-1&2&-1&-1&0 \\
        5&3&-5&-1&3&-1&-2&0 \\ 4&3&-4&-1&2&-1&-1&0 \\
        3&2&-3&0&1&-1&-1&0 \\ 2&2&-2&0&0&0&-1&0 \\ 1&1&-1&0&0&0&0&-1 \end{array}\right)$,
  $\left(\begin{array}{cccccccc} 0&-2&1&0&0&1&-1&0 \\ -1&-2&1&0&0&2&-2&0 \\ 0&-3&1&0&0&3&-3&0 \\
        -1&-5&1&1&0&4&-4&0 \\ -1&-4&1&0&1&3&-3&0 \\ -1&-3&1&0&0&3&-2&0 \\
        0&-2&0&0&0&2&-1&0 \\ 0&-1&0&0&0&1&-1&1 \end{array}\right)$,
  $\left(\begin{array}{cccccccc} 1&0&-2&0&1&1&-1&0 \\ 3&0&-3&0&1&1&-1&0 \\ 3&0&-4&0&2&1&-1&0 \\
        5&0&-5&-1&3&2&-2&0 \\ 4&0&-4&-1&2&2&-1&0 \\ 3&-1&-3&0&1&2&-1&0 \\
        2&0&-2&0&0&2&-1&0 \\ 1&0&-1&0&0&1&-1&1 \end{array}\right)$,
  $\left(\begin{array}{cccccccc} 0&-2&1&0&0&1&-1&0 \\ -1&-2&1&0&0&2&-2&0 \\ 0&-3&1&0&0&3&-3&0 \\
        -1&-5&1&1&0&4&-4&0 \\ -1&-4&1&0&1&3&-3&0 \\ -1&-3&1&0&0&3&-2&0 \\
        0&-2&0&0&0&2&-1&0 \\ 0&-1&0&0&0&1&0&-1 \end{array}\right)$,
  $\left(\begin{array}{cccccccc} 1&0&-2&0&1&1&-1&0 \\ 3&0&-3&0&1&1&-1&0 \\ 3&0&-4&0&2&1&-1&0 \\
        5&0&-5&-1&3&2&-2&0 \\ 4&0&-4&-1&2&2&-1&0 \\ 3&-1&-3&0&1&2&-1&0 \\
        2&0&-2&0&0&2&-1&0 \\ 1&0&-1&0&0&1&0&-1 \end{array}\right)$,
  $\left(\begin{array}{cccccccc} 0&-1&1&0&0&0&-1&0 \\ -1&-1&1&0&0&1&-2&0 \\ 0&-1&1&0&0&1&-3&0 \\
        -1&-2&1&1&0&1&-4&0 \\ -1&-1&1&0&1&0&-3&0 \\ -1&0&1&0&0&0&-2&0 \\
        0&0&0&0&0&0&-1&0 \\ 0&0&0&0&0&0&0&-1 \end{array}\right)$,
  $\left(\begin{array}{cccccccc} 1&1&-2&0&1&0&-1&0 \\ 3&1&-3&0&1&0&-1&0 \\ 3&2&-4&0&2&-1&-1&0 \\
        5&3&-5&-1&3&-1&-2&0 \\ 4&3&-4&-1&2&-1&-1&0 \\
        3&2&-3&0&1&-1&-1&0 \\ 2&2&-2&0&0&0&-1&0 \\ 1&1&-1&0&0&0&0&-1 \end{array}\right)$,
  $\left(\begin{array}{cccccccc} 1&0&-1&0&0&2&-2&0 \\ 2&0&-3&1&0&2&-2&0 \\ 2&0&-3&1&0&3&-3&0 \\
        4&-1&-5&2&0&4&-4&0 \\ 3&-1&-4&1&1&3&-3&0 \\ 2&-1&-3&1&0&3&-2&0 \\
        2&0&-2&0&0&2&-1&0 \\ 1&0&-1&0&0&1&-1&1 \end{array}\right)$,
  $\left(\begin{array}{cccccccc} -2&0&2&-1&0&0&1&0 \\ -3&-1&2&0&0&-1&2&0 \\ -4&-1&3&-1&0&0&2&0 \\
        -5&-2&4&-1&0&-1&4&0 \\ -4&-2&3&0&-1&0&3&0 \\ -3&-2&2&0&0&0&2&0 \\
        -2&-2&2&0&0&0&1&0 \\ -1&-1&1&0&0&0&0&1 \end{array}\right)$,
  $\left(\begin{array}{cccccccc} 1&0&-1&0&1&-1&1&0 \\ 1&1&0&-1&1&-1&1&0 \\ 1&0&0&-1&2&-2&2&0 \\
        1&1&0&-2&3&-2&2&0 \\ 1&1&0&-2&2&-1&2&0 \\ 1&0&0&-1&1&0&1&0 \\
        0&0&0&0&0&0&1&0 \\ 0&0&0&0&0&0&0&1 \end{array}\right)$,
  $\left(\begin{array}{cccccccc} 0&-2&0&1&-1&1&0&0 \\ 0&-2&1&0&-1&2&-1&0 \\ 1&-3&0&1&-2&3&-1&0 \\
        0&-4&1&1&-3&5&-2&0 \\ 0&-4&1&1&-2&4&-2&0 \\ 0&-3&1&0&-1&3&-1&0 \\
        0&-2&0&0&0&2&-1&0 \\ 0&-1&0&0&0&1&-1&1 \end{array}\right)$,
  $\left(\begin{array}{cccccccc} 1&0&-1&0&0&2&-2&0 \\ 2&0&-3&1&0&2&-2&0 \\ 2&0&-3&1&0&3&-3&0 \\
        4&-1&-5&2&0&4&-4&0 \\ 3&-1&-4&1&1&3&-3&0 \\ 2&-1&-3&1&0&3&-2&0 \\
        2&0&-2&0&0&2&-1&0 \\ 1&0&-1&0&0&1&0&-1 \end{array}\right)$,
  $\left(\begin{array}{cccccccc} -2&0&2&-1&0&0&1&0 \\ -3&-1&2&0&0&-1&2&0 \\ -4&-1&3&-1&0&0&2&0 \\
        -5&-2&4&-1&0&-1&4&0 \\ -4&-2&3&0&-1&0&3&0 \\ -3&-2&2&0&0&0&2&0 \\
        -2&-2&2&0&0&0&1&0 \\ -1&-1&1&0&0&0&1&-1 \end{array}\right)$,
  $\left(\begin{array}{cccccccc} 1&0&-1&0&1&-1&1&0 \\ 1&1&0&-1&1&-1&1&0 \\ 1&0&0&-1&2&-2&2&0 \\
        1&1&0&-2&3&-2&2&0 \\ 1&1&0&-2&2&-1&2&0 \\ 1&0&0&-1&1&0&1&0 \\
        0&0&0&0&0&0&1&0 \\ 0&0&0&0&0&0&1&-1 \end{array}\right)$,
  $\left(\begin{array}{cccccccc} 0&-2&0&1&-1&1&0&0 \\ 0&-2&1&0&-1&2&-1&0 \\ 1&-3&0&1&-2&3&-1&0 \\
        0&-4&1&1&-3&5&-2&0 \\ 0&-4&1&1&-2&4&-2&0 \\ 0&-3&1&0&-1&3&-1&0 \\
        0&-2&0&0&0&2&-1&0 \\ 0&-1&0&0&0&1&0&-1 \end{array}\right)$,
  $\left(\begin{array}{cccccccc} 1&1&-1&0&0&1&-2&0 \\ 2&1&-3&1&0&1&-2&0 \\ 2&2&-3&1&0&1&-3&0 \\
        4&2&-5&2&0&1&-4&0 \\ 3&2&-4&1&1&0&-3&0 \\ 2&2&-3&1&0&0&-2&0 \\
        2&2&-2&0&0&0&-1&0 \\ 1&1&-1&0&0&0&0&-1 \end{array}\right)$,
  $\left(\begin{array}{cccccccc} -2&1&2&-1&0&-1&1&0 \\ -3&0&2&0&0&-2&2&0 \\ -4&1&3&-1&0&-2&2&0 \\
        -5&1&4&-1&0&-4&4&0 \\ -4&1&3&0&-1&-3&3&0 \\ -3&1&2&0&0&-3&2&0 \\
        -2&0&2&0&0&-2&1&0 \\ -1&0&1&0&0&-1&1&-1 \end{array}\right)$,
  $\left(\begin{array}{cccccccc} 1&1&-1&0&1&-2&1&0 \\ 1&2&0&-1&1&-2&1&0 \\ 1&2&0&-1&2&-4&2&0 \\
        1&4&0&-2&3&-5&2&0 \\ 1&4&0&-2&2&-4&2&0 \\ 1&3&0&-1&1&-3&1&0 \\
        0&2&0&0&0&-2&1&0 \\ 0&1&0&0&0&-1&1&-1 \end{array}\right)$,
  $\left(\begin{array}{cccccccc} 0&-1&0&1&-1&0&0&0 \\ 0&-1&1&0&-1&1&-1&0 \\ 1&-1&0&1&-2&1&-1&0 \\
        0&-1&1&1&-3&2&-2&0 \\ 0&-1&1&1&-2&1&-2&0 \\ 0&0&1&0&-1&0&-1&0 \\
        0&0&0&0&0&0&-1&0 \\ 0&0&0&0&0&0&0&-1 \end{array}\right)$,
  $\left(\begin{array}{cccccccc} 2&-1&-2&1&0&1&-1&0 \\ 2&0&-3&1&0&2&-2&0 \\ 3&-1&-4&2&0&2&-2&0 \\
        4&-1&-5&2&0&4&-4&0 \\ 3&-1&-4&1&1&3&-3&0 \\ 2&-1&-3&1&0&3&-2&0 \\
        2&0&-2&0&0&2&-1&0 \\ 1&0&-1&0&0&1&-1&1 \end{array}\right)$,
  $\left(\begin{array}{cccccccc} -1&-1&1&0&0&-1&2&0 \\ -3&-1&2&0&0&-1&2&0 \\ -3&-2&2&0&0&-1&3&0 \\
        -5&-2&4&-1&0&-1&4&0 \\ -4&-2&3&0&-1&0&3&0 \\ -3&-2&2&0&0&0&2&0 \\
        -2&-2&2&0&0&0&1&0 \\ -1&-1&1&0&0&0&0&1 \end{array}\right)$,
  $\left(\begin{array}{cccccccc} 0&1&0&-1&1&0&0&0 \\ 1&1&0&-1&1&-1&1&0 \\ 0&1&1&-2&2&-1&1&0 \\
        1&1&0&-2&3&-2&2&0 \\ 1&1&0&-2&2&-1&2&0 \\ 1&0&0&-1&1&0&1&0 \\
        0&0&0&0&0&0&1&0 \\ 0&0&0&0&0&0&0&1 \end{array}\right)$,
  $\left(\begin{array}{cccccccc} -1&-1&1&0&-1&2&-1&0 \\ 0&-2&1&0&-1&2&-1&0 \\ 0&-2&1&0&-2&4&-2&0 \\
        0&-4&1&1&-3&5&-2&0 \\ 0&-4&1&1&-2&4&-2&0 \\ 0&-3&1&0&-1&3&-1&0 \\
        0&-2&0&0&0&2&-1&0 \\ 0&-1&0&0&0&1&-1&1 \end{array}\right)$,
  $\left(\begin{array}{cccccccc} 2&-1&-2&1&0&1&-1&0 \\ 2&0&-3&1&0&2&-2&0 \\ 3&-1&-4&2&0&2&-2&0 \\
        4&-1&-5&2&0&4&-4&0 \\ 3&-1&-4&1&1&3&-3&0 \\ 2&-1&-3&1&0&3&-2&0 \\
        2&0&-2&0&0&2&-1&0 \\ 1&0&-1&0&0&1&0&-1 \end{array}\right)$,
  $\left(\begin{array}{cccccccc} -1&-1&1&0&0&-1&2&0 \\ -3&-1&2&0&0&-1&2&0 \\ -3&-2&2&0&0&-1&3&0 \\
        -5&-2&4&-1&0&-1&4&0 \\ -4&-2&3&0&-1&0&3&0 \\ -3&-2&2&0&0&0&2&0 \\
        -2&-2&2&0&0&0&1&0 \\ -1&-1&1&0&0&0&1&-1 \end{array}\right)$,
  $\left(\begin{array}{cccccccc} 0&1&0&-1&1&0&0&0 \\ 1&1&0&-1&1&-1&1&0 \\ 0&1&1&-2&2&-1&1&0 \\
        1&1&0&-2&3&-2&2&0 \\ 1&1&0&-2&2&-1&2&0 \\ 1&0&0&-1&1&0&1&0 \\
        0&0&0&0&0&0&1&0 \\ 0&0&0&0&0&0&1&-1 \end{array}\right)$,
  $\left(\begin{array}{cccccccc} -1&-1&1&0&-1&2&-1&0 \\ 0&-2&1&0&-1&2&-1&0 \\ 0&-2&1&0&-2&4&-2&0 \\
        0&-4&1&1&-3&5&-2&0 \\ 0&-4&1&1&-2&4&-2&0 \\ 0&-3&1&0&-1&3&-1&0 \\
        0&-2&0&0&0&2&-1&0 \\ 0&-1&0&0&0&1&0&-1 \end{array}\right)$,
  $\left(\begin{array}{cccccccc} 2&0&-2&1&0&0&-1&0 \\ 2&1&-3&1&0&1&-2&0 \\ 3&1&-4&2&0&0&-2&0 \\
        4&2&-5&2&0&1&-4&0 \\ 3&2&-4&1&1&0&-3&0 \\ 2&2&-3&1&0&0&-2&0 \\
        2&2&-2&0&0&0&-1&0 \\ 1&1&-1&0&0&0&0&-1 \end{array}\right)$,
  $\left(\begin{array}{cccccccc} -1&0&1&0&0&-2&2&0 \\ -3&0&2&0&0&-2&2&0 \\ -3&0&2&0&0&-3&3&0 \\
        -5&1&4&-1&0&-4&4&0 \\ -4&1&3&0&-1&-3&3&0 \\ -3&1&2&0&0&-3&2&0 \\
        -2&0&2&0&0&-2&1&0 \\ -1&0&1&0&0&-1&1&-1 \end{array}\right)$,
  $\left(\begin{array}{cccccccc} 0&2&0&-1&1&-1&0&0 \\ 1&2&0&-1&1&-2&1&0 \\ 0&3&1&-2&2&-3&1&0 \\
        1&4&0&-2&3&-5&2&0 \\ 1&4&0&-2&2&-4&2&0 \\ 1&3&0&-1&1&-3&1&0 \\
        0&2&0&0&0&-2&1&0 \\ 0&1&0&0&0&-1&1&-1 \end{array}\right)$,
  $\left(\begin{array}{cccccccc} -1&0&1&0&-1&1&-1&0 \\ 0&-1&1&0&-1&1&-1&0 \\ 0&0&1&0&-2&2&-2&0 \\
        0&-1&1&1&-3&2&-2&0 \\ 0&-1&1&1&-2&1&-2&0 \\ 0&0&1&0&-1&0&-1&0 \\
        0&0&0&0&0&0&-1&0 \\ 0&0&0&0&0&0&0&-1 \end{array}\right)$,
  $\left(\begin{array}{cccccccc} 1&1&0&-1&1&-1&0&1 \\ 0&2&0&-1&1&-1&0&1 \\ 0&2&1&-2&2&-2&0&2 \\
        0&3&0&-2&3&-3&0&3 \\ 0&3&0&-2&2&-2&0&3 \\ 0&2&0&-1&1&-1&-1&3 \\
        0&1&0&0&0&-1&0&2 \\ 0&0&0&0&0&0&0&1 \end{array}\right)$,
  $\left(\begin{array}{cccccccc} 0&-1&1&0&-1&1&-1&1 \\ -1&-1&1&0&-1&2&-2&1 \\ 0&-1&1&0&-2&3&-3&2 \\
        -1&-2&1&1&-3&4&-4&3 \\ -1&-2&1&1&-2&3&-4&3 \\
        -1&-1&1&0&-1&2&-3&3 \\ 0&-1&0&0&0&1&-2&2 \\ 0&-1&0&0&0&1&-1&1\end{array}\right)$,
  $\left(\begin{array}{cccccccc} -2&0&1&0&0&-1&1&1 \\ -2&-1&2&0&0&-2&2&1 \\ -3&-1&2&0&0&-2&2&2 \\
        -4&-1&4&-1&0&-3&3&3 \\ -3&-1&3&0&-1&-2&2&3 \\
        -2&-1&2&0&0&-2&1&3 \\ -2&-1&2&0&0&-1&0&2 \\ -1&-1&1&0&0&0&0&1 \end{array}\right)$,
  $\left(\begin{array}{cccccccc} 1&0&-2&1&0&1&-2&1 \\ 3&0&-3&1&0&1&-2&1 \\ 3&0&-4&2&0&1&-3&2 \\
        5&0&-5&2&0&2&-5&3 \\ 4&0&-4&1&1&1&-4&3 \\ 3&0&-3&1&0&1&-3&3 \\
        2&1&-2&0&0&1&-2&2 \\ 1&0&-1&0&0&1&-1&1 \end{array}\right)$,
  $\left(\begin{array}{cccccccc} 1&1&0&-1&1&-1&1&-1 \\ 0&2&0&-1&1&-1&1&-1 \\ 0&2&1&-2&2&-2&2&-2 \\
        0&3&0&-2&3&-3&3&-3 \\ 0&3&0&-2&2&-2&3&-3 \\ 0&2&0&-1&1&-1&2&-3 \\
        0&1&0&0&0&-1&2&-2 \\ 0&0&0&0&0&0&1&-1 \end{array}\right)$,
  $\left(\begin{array}{cccccccc} 0&-1&1&0&-1&1&0&-1 \\ -1&-1&1&0&-1&2&-1&-1 \\
        0&-1&1&0&-2&3&-1&-2 \\ -1&-2&1&1&-3&4&-1&-3 \\
        -1&-2&1&1&-2&3&-1&-3 \\ -1&-1&1&0&-1&2&0&-3 \\
        0&-1&0&0&0&1&0&-2 \\ 0&-1&0&0&0&1&0&-1 \end{array}\right)$,
  $\left(\begin{array}{cccccccc} -2&0&1&0&0&-1&2&-1 \\ -2&-1&2&0&0&-2&3&-1 \\ -3&-1&2&0&0&-2&4&-2 \\
        -4&-1&4&-1&0&-3&6&-3 \\ -3&-1&3&0&-1&-2&5&-3 \\
        -2&-1&2&0&0&-2&4&-3 \\ -2&-1&2&0&0&-1&2&-2 \\
        -1&-1&1&0&0&0&1&-1 \end{array}\right)$,
  $\left(\begin{array}{cccccccc} 1&0&-2&1&0&1&-1&-1 \\ 3&0&-3&1&0&1&-1&-1 \\ 3&0&-4&2&0&1&-1&-2 \\
        5&0&-5&2&0&2&-2&-3 \\ 4&0&-4&1&1&1&-1&-3 \\ 3&0&-3&1&0&1&0&-3 \\
        2&1&-2&0&0&1&0&-2 \\ 1&0&-1&0&0&1&0&-1 \end{array}\right)$,
  $\left(\begin{array}{cccccccc} 1&1&0&-1&1&-1&0&1 \\ 1&1&-1&0&1&-2&1&1 \\ 1&1&0&-1&2&-3&1&2 \\
        1&2&-1&-1&3&-4&1&3 \\ 1&2&-1&-1&2&-3&1&3 \\ 1&1&-1&0&1&-2&0&3 \\
        0&1&0&0&0&-1&0&2 \\ 0&0&0&0&0&0&0&1 \end{array}\right)$,
  $\left(\begin{array}{cccccccc} 0&-1&1&0&-1&1&-1&1 \\ 0&-2&0&1&-1&1&-1&1 \\ 1&-2&0&1&-2&2&-2&2 \\
        0&-3&0&2&-3&3&-3&3 \\ 0&-3&0&2&-2&2&-3&3 \\ 0&-2&0&1&-1&1&-2&3 \\
        0&-1&0&0&0&1&-2&2 \\ 0&-1&0&0&0&1&-1&1 \end{array}\right)$,
  $\left(\begin{array}{cccccccc} -2&0&1&0&0&-1&1&1 \\ -3&0&3&-1&0&-1&1&1 \\ -4&0&3&-1&0&-1&1&2 \\
        -5&0&5&-2&0&-2&2&3 \\ -4&0&4&-1&-1&-1&1&3 \\
        -3&0&3&-1&0&-1&0&3 \\ -2&-1&2&0&0&-1&0&2 \\ -1&-1&1&0&0&0&0&1 \end{array}\right)$,
  $\left(\begin{array}{cccccccc} 1&0&-2&1&0&1&-2&1 \\ 2&1&-2&0&0&2&-3&1 \\ 2&1&-3&1&0&2&-4&2 \\
        4&1&-4&1&0&3&-6&3 \\ 3&1&-3&0&1&2&-5&3 \\ 2&1&-2&0&0&2&-4&3 \\
        2&1&-2&0&0&1&-2&2 \\ 1&0&-1&0&0&1&-1&1 \end{array}\right)$,
  $\left(\begin{array}{cccccccc} 1&1&0&-1&1&-1&1&-1 \\ 1&1&-1&0&1&-2&2&-1 \\ 1&1&0&-1&2&-3&3&-2 \\
        1&2&-1&-1&3&-4&4&-3 \\ 1&2&-1&-1&2&-3&4&-3 \\
        1&1&-1&0&1&-2&3&-3 \\ 0&1&0&0&0&-1&2&-2 \\ 0&0&0&0&0&0&1&-1 \end{array}\right)$,
  $\left(\begin{array}{cccccccc} 0&-1&1&0&-1&1&0&-1 \\ 0&-2&0&1&-1&1&0&-1 \\
        1&-2&0&1&-2&2&0&-2 \\ 0&-3&0&2&-3&3&0&-3 \\ 0&-3&0&2&-2&2&0&-3 \\
        0&-2&0&1&-1&1&1&-3 \\ 0&-1&0&0&0&1&0&-2 \\ 0&-1&0&0&0&1&0&-1 \end{array}\right)$,
  $\left(\begin{array}{cccccccc} -2&0&1&0&0&-1&2&-1 \\ -3&0&3&-1&0&-1&2&-1 \\
        -4&0&3&-1&0&-1&3&-2 \\ -5&0&5&-2&0&-2&5&-3 \\
        -4&0&4&-1&-1&-1&4&-3 \\ -3&0&3&-1&0&-1&3&-3 \\
        -2&-1&2&0&0&-1&2&-2 \\ -1&-1&1&0&0&0&1&-1 \end{array}\right)$,
  $\left(\begin{array}{cccccccc} 1&0&-2&1&0&1&-1&-1 \\ 2&1&-2&0&0&2&-2&-1 \\ 2&1&-3&1&0&2&-2&-2 \\
        4&1&-4&1&0&3&-3&-3 \\ 3&1&-3&0&1&2&-2&-3 \\ 2&1&-2&0&0&2&-1&-3 \\
        2&1&-2&0&0&1&0&-2 \\ 1&0&-1&0&0&1&0&-1 \end{array}\right)$,
  $\left(\begin{array}{cccccccc} 0&1&-1&0&1&-1&0&1 \\ 1&1&-1&0&1&-2&1&1 \\ 0&1&-1&0&2&-3&1&2 \\
        1&2&-1&-1&3&-4&1&3 \\ 1&2&-1&-1&2&-3&1&3 \\ 1&1&-1&0&1&-2&0&3 \\
        0&1&0&0&0&-1&0&2 \\ 0&0&0&0&0&0&0&1 \end{array}\right)$,
  $\left(\begin{array}{cccccccc} -1&-1&0&1&-1&1&-1&1 \\ 0&-2&0&1&-1&1&-1&1 \\ 0&-2&-1&2&-2&2&-2&2 \\
        0&-3&0&2&-3&3&-3&3 \\ 0&-3&0&2&-2&2&-3&3 \\ 0&-2&0&1&-1&1&-2&3 \\
        0&-1&0&0&0&1&-2&2 \\ 0&-1&0&0&0&1&-1&1 \end{array}\right)$,
  $\left(\begin{array}{cccccccc} -1&0&2&-1&0&-1&1&1 \\ -3&0&3&-1&0&-1&1&1 \\ -3&0&4&-2&0&-1&1&2 \\
        -5&0&5&-2&0&-2&2&3 \\ -4&0&4&-1&-1&-1&1&3 \\
        -3&0&3&-1&0&-1&0&3 \\ -2&-1&2&0&0&-1&0&2 \\ -1&-1&1&0&0&0&0&1 \end{array}\right)$,
  $\left(\begin{array}{cccccccc} 2&0&-1&0&0&1&-2&1 \\ 2&1&-2&0&0&2&-3&1 \\ 3&1&-2&0&0&2&-4&2 \\
        4&1&-4&1&0&3&-6&3 \\ 3&1&-3&0&1&2&-5&3 \\ 2&1&-2&0&0&2&-4&3 \\
        2&1&-2&0&0&1&-2&2 \\ 1&0&-1&0&0&1&-1&1 \end{array}\right)$,
  $\left(\begin{array}{cccccccc} 0&1&-1&0&1&-1&1&-1 \\ 1&1&-1&0&1&-2&2&-1 \\ 0&1&-1&0&2&-3&3&-2 \\
        1&2&-1&-1&3&-4&4&-3 \\ 1&2&-1&-1&2&-3&4&-3 \\
        1&1&-1&0&1&-2&3&-3 \\ 0&1&0&0&0&-1&2&-2 \\ 0&0&0&0&0&0&1&-1 \end{array}\right)$,
   $\left(\begin{array}{cccccccc} -1&-1&0&1&-1&1&0&-1 \\ 0&-2&0&1&-1&1&0&-1 \\
        0&-2&-1&2&-2&2&0&-2 \\ 0&-3&0&2&-3&3&0&-3 \\
        0&-3&0&2&-2&2&0&-3 \\ 0&-2&0&1&-1&1&1&-3 \\ 0&-1&0&0&0&1&0&-2 \\
        0&-1&0&0&0&1&0&-1 \end{array}\right)$,
  $\left(\begin{array}{cccccccc} -1&0&2&-1&0&-1&2&-1 \\ -3&0&3&-1&0&-1&2&-1 \\
        -3&0&4&-2&0&-1&3&-2 \\ -5&0&5&-2&0&-2&5&-3 \\
        -4&0&4&-1&-1&-1&4&-3 \\ -3&0&3&-1&0&-1&3&-3 \\
        -2&-1&2&0&0&-1&2&-2 \\ -1&-1&1&0&0&0&1&-1 \end{array}\right)$,
  $\left(\begin{array}{cccccccc} 2&0&-1&0&0&1&-1&-1 \\ 2&1&-2&0&0&2&-2&-1 \\ 3&1&-2&0&0&2&-2&-2 \\
        4&1&-4&1&0&3&-3&-3 \\ 3&1&-3&0&1&2&-2&-3 \\ 2&1&-2&0&0&2&-1&-3 \\
        2&1&-2&0&0&1&0&-2 \\ 1&0&-1&0&0&1&0&-1 \end{array}\right)$,
  $\left(\begin{array}{cccccccc} 1&1&-1&-1&1&1&-2&1 \\ 2&1&-3&0&1&1&-2&1 \\ 2&2&-3&-1&2&1&-3&2 \\
        4&2&-5&-1&3&1&-4&3 \\ 3&2&-4&-1&2&1&-3&3 \\ 2&1&-3&0&1&1&-3&3 \\
        2&1&-2&0&0&1&-2&2 \\ 1&0&-1&0&0&1&-1&1 \end{array}\right)$,
  $\left(\begin{array}{cccccccc} 1&1&-1&0&0&-1&1&1 \\ 1&2&0&-1&0&-1&1&1 \\ 1&2&0&-1&0&-2&2&2 \\
        1&4&0&-2&0&-2&2&3 \\ 1&3&0&-1&-1&-1&1&3 \\ 1&2&0&-1&0&-1&0&3 \\
        0&1&0&0&0&-1&0&2 \\ 0&0&0&0&0&0&0&1 \end{array}\right)$,
  $\left(\begin{array}{cccccccc} 1&1&-1&-1&1&1&-1&-1 \\ 2&1&-3&0&1&1&-1&-1 \\ 2&2&-3&-1&2&1&-1&-2 \\
        4&2&-5&-1&3&1&-1&-3 \\ 3&2&-4&-1&2&1&0&-3 \\ 2&1&-3&0&1&1&0&-3 \\
        2&1&-2&0&0&1&0&-2 \\ 1&0&-1&0&0&1&0&-1 \end{array}\right)$,
  $\left(\begin{array}{cccccccc} 1&1&-1&0&0&-1&2&-1 \\ 1&2&0&-1&0&-1&2&-1 \\ 1&2&0&-1&0&-2&4&-2 \\
        1&4&0&-2&0&-2&5&-3 \\ 1&3&0&-1&-1&-1&4&-3 \\
        1&2&0&-1&0&-1&3&-3 \\ 0&1&0&0&0&-1&2&-2 \\ 0&0&0&0&0&0&1&-1 \end{array}\right)$,
  $\left(\begin{array}{cccccccc} 2&0&-2&0&1&0&-1&1 \\ 2&1&-3&0&1&1&-2&1 \\ 3&1&-4&0&2&0&-2&2 \\
        4&2&-5&-1&3&1&-4&3 \\ 3&2&-4&-1&2&1&-3&3 \\ 2&1&-3&0&1&1&-3&3 \\
        2&1&-2&0&0&1&-2&2 \\ 1&0&-1&0&0&1&-1&1 \end{array}\right)$,
  $\left(\begin{array}{cccccccc} 0&2&0&-1&0&0&0&1 \\ 1&2&0&-1&0&-1&1&1 \\ 0&3&1&-2&0&-1&1&2 \\
        1&4&0&-2&0&-2&2&3 \\ 1&3&0&-1&-1&-1&1&3 \\ 1&2&0&-1&0&-1&0&3 \\
        0&1&0&0&0&-1&0&2 \\ 0&0&0&0&0&0&0&1 \end{array}\right)$,
  $\left(\begin{array}{cccccccc} 2&0&-2&0&1&0&0&-1 \\ 2&1&-3&0&1&1&-1&-1 \\ 3&1&-4&0&2&0&0&-2 \\
        4&2&-5&-1&3&1&-1&-3 \\ 3&2&-4&-1&2&1&0&-3 \\ 2&1&-3&0&1&1&0&-3 \\
        2&1&-2&0&0&1&0&-2 \\ 1&0&-1&0&0&1&0&-1 \end{array}\right)$,
  $\left(\begin{array}{cccccccc} 0&2&0&-1&0&0&1&-1 \\ 1&2&0&-1&0&-1&2&-1 \\ 0&3&1&-2&0&-1&3&-2 \\
        1&4&0&-2&0&-2&5&-3 \\ 1&3&0&-1&-1&-1&4&-3 \\
        1&2&0&-1&0&-1&3&-3 \\ 0&1&0&0&0&-1&2&-2 \\ 0&0&0&0&0&0&1&-1 \end{array}\right)$,
  $\left(\begin{array}{cccccccc} -1&-1&1&0&1&-2&1&1 \\ -2&-1&3&-1&1&-2&1&1 \\
        -2&-2&3&-1&2&-3&1&2 \\ -4&-2&5&-2&3&-4&1&3 \\
        -3&-1&4&-2&2&-3&1&3 \\ -2&-1&3&-1&1&-2&0&3 \\
        -2&-1&2&0&0&-1&0&2 \\ -1&-1&1&0&0&0&0&1 \end{array}\right)$,
  $\left(\begin{array}{cccccccc} 2&0&-2&1&-1&1&-1&1 \\ 3&1&-2&0&-1&2&-2&1 \\ 4&1&-3&1&-2&2&-2&2 \\
        5&2&-4&1&-3&4&-4&3 \\ 4&1&-3&1&-2&3&-4&3 \\ 3&1&-2&0&-1&2&-3&3 \\
        2&1&-2&0&0&1&-2&2 \\ 1&0&-1&0&0&1&-1&1 \end{array}\right)$,
  $\left(\begin{array}{cccccccc} -1&0&1&-1&0&1&-1&1 \\ -1&-1&0&0&0&1&-1&1 \\ -1&0&0&-1&0&2&-2&2 \\
        -1&-1&0&-1&0&2&-2&3 \\ -1&-1&0&0&-1&2&-2&3 \\
        -1&-1&0&0&0&1&-2&3 \\ 0&-1&0&0&0&1&-2&2 \\ 0&-1&0&0&0&1&-1&1 \end{array}\right)$,
  $\left(\begin{array}{cccccccc} 0&1&0&0&0&0&-1&1 \\ 0&1&-1&1&0&-1&0&1 \\ -1&1&0&1&0&-1&-1&2 \\
        0&1&-1&2&0&-2&-1&3 \\ 0&1&-1&1&1&-2&-1&3 \\ 0&1&-1&1&0&-1&-1&3 \\
        0&1&0&0&0&-1&0&2 \\ 0&0&0&0&0&0&0&1 \end{array}\right)$,
  $\left(\begin{array}{cccccccc} -1&-1&1&0&1&-2&2&-1 \\ -2&-1&3&-1&1&-2&2&-1 \\
        -2&-2&3&-1&2&-3&3&-2 \\ -4&-2&5&-2&3&-4&4&-3 \\
        -3&-1&4&-2&2&-3&4&-3 \\ -2&-1&3&-1&1&-2&3&-3 \\
        -2&-1&2&0&0&-1&2&-2 \\ -1&-1&1&0&0&0&1&-1 \end{array}\right)$,
  $\left(\begin{array}{cccccccc} 2&0&-2&1&-1&1&0&-1 \\ 3&1&-2&0&-1&2&-1&-1 \\ 4&1&-3&1&-2&2&0&-2 \\
        5&2&-4&1&-3&4&-1&-3 \\ 4&1&-3&1&-2&3&-1&-3 \\
        3&1&-2&0&-1&2&0&-3 \\ 2&1&-2&0&0&1&0&-2 \\ 1&0&-1&0&0&1&0&-1 \end{array}\right)$,
  $\left(\begin{array}{cccccccc} -1&0&1&-1&0&1&0&-1 \\ -1&-1&0&0&0&1&0&-1 \\
        -1&0&0&-1&0&2&0&-2 \\ -1&-1&0&-1&0&2&1&-3 \\
        -1&-1&0&0&-1&2&1&-3 \\ -1&-1&0&0&0&1&1&-3 \\ 0&-1&0&0&0&1&0&-2 \\
        0&-1&0&0&0&1&0&-1 \end{array}\right)$,
  $\left(\begin{array}{cccccccc} 0&1&0&0&0&0&0&-1 \\ 0&1&-1&1&0&-1&1&-1 \\ -1&1&0&1&0&-1&1&-2 \\
        0&1&-1&2&0&-2&2&-3 \\ 0&1&-1&1&1&-2&2&-3 \\ 0&1&-1&1&0&-1&2&-3 \\
        0&1&0&0&0&-1&2&-2 \\ 0&0&0&0&0&0&1&-1 \end{array}\right)$,
  $\left(\begin{array}{cccccccc} -2&0&2&-1&1&-1&0&1 \\ -2&-1&3&-1&1&-2&1&1 \\ -3&-1&4&-2&2&-2&0&2 \\
        -4&-2&5&-2&3&-4&1&3 \\ -3&-1&4&-2&2&-3&1&3 \\
        -2&-1&3&-1&1&-2&0&3 \\ -2&-1&2&0&0&-1&0&2 \\ -1&-1&1&0&0&0&0&1 \end{array}\right)$,
  $\left(\begin{array}{cccccccc} 1&1&-1&0&-1&2&-2&1 \\ 3&1&-2&0&-1&2&-2&1 \\
        3&2&-2&0&-2&3&-3&2 \\ 5&2&-4&1&-3&4&-4&3 \\ 4&1&-3&1&-2&3&-4&3 \\
        3&1&-2&0&-1&2&-3&3 \\ 2&1&-2&0&0&1&-2&2 \\ 1&0&-1&0&0&1&-1&1 \end{array}\right)$,
  $\left(\begin{array}{cccccccc} 0&-1&0&0&0&0&0&1 \\ -1&-1&0&0&0&1&-1&1 \\ 0&-1&-1&0&0&1&-1&2 \\
        -1&-1&0&-1&0&2&-2&3 \\ -1&-1&0&0&-1&2&-2&3 \\
        -1&-1&0&0&0&1&-2&3 \\ 0&-1&0&0&0&1&-2&2 \\ 0&-1&0&0&0&1&-1&1 \end{array}\right)$,
  $\left(\begin{array}{cccccccc} 1&0&-1&1&0&-1&0&1 \\ 0&1&-1&1&0&-1&0&1 \\ 0&0&-1&2&0&-2&0&2 \\
        0&1&-1&2&0&-2&-1&3 \\ 0&1&-1&1&1&-2&-1&3 \\ 0&1&-1&1&0&-1&-1&3 \\
        0&1&0&0&0&-1&0&2 \\ 0&0&0&0&0&0&0&1 \end{array}\right)$,
  $\left(\begin{array}{cccccccc} -2&0&2&-1&1&-1&1&-1 \\ -2&-1&3&-1&1&-2&2&-1 \\
        -3&-1&4&-2&2&-2&2&-2 \\ -4&-2&5&-2&3&-4&4&-3 \\
        -3&-1&4&-2&2&-3&4&-3 \\ -2&-1&3&-1&1&-2&3&-3 \\
        -2&-1&2&0&0&-1&2&-2 \\ -1&-1&1&0&0&0&1&-1 \end{array}\right)$,
  $\left(\begin{array}{cccccccc} 1&1&-1&0&-1&2&-1&-1 \\ 3&1&-2&0&-1&2&-1&-1 \\
        3&2&-2&0&-2&3&-1&-2 \\ 5&2&-4&1&-3&4&-1&-3 \\
        4&1&-3&1&-2&3&-1&-3 \\ 3&1&-2&0&-1&2&0&-3 \\ 2&1&-2&0&0&1&0&-2 \\
        1&0&-1&0&0&1&0&-1 \end{array}\right)$,
  $\left(\begin{array}{cccccccc} 0&-1&0&0&0&0&1&-1 \\ -1&-1&0&0&0&1&0&-1 \\ 0&-1&-1&0&0&1&1&-2 \\
        -1&-1&0&-1&0&2&1&-3 \\ -1&-1&0&0&-1&2&1&-3 \\
        -1&-1&0&0&0&1&1&-3 \\ 0&-1&0&0&0&1&0&-2 \\ 0&-1&0&0&0&1&0&-1 \end{array}\right)$,
  $\left(\begin{array}{cccccccc} 1&0&-1&1&0&-1&1&-1 \\ 0&1&-1&1&0&-1&1&-1 \\
        0&0&-1&2&0&-2&2&-2 \\ 0&1&-1&2&0&-2&2&-3 \\ 0&1&-1&1&1&-2&2&-3 \\
        0&1&-1&1&0&-1&2&-3 \\ 0&1&0&0&0&-1&2&-2 \\ 0&0&0&0&0&0&1&-1 \end{array}\right)$,
  $\left(\begin{array}{cccccccc} -1&0&0&2&-2&0&0&0 \\ -1&0&1&2&-3&1&-1&0 \\ -2&0&1&3&-3&0&-1&0 \\
        -2&0&2&4&-5&0&-1&0 \\ -2&0&2&3&-4&0&-1&0 \\ -1&1&1&2&-3&0&-1&0 \\
        0&0&0&2&-2&0&-1&0 \\ 0&0&0&1&-1&0&-1&1 \end{array}\right)$,
  $\left(\begin{array}{cccccccc} 0&-1&1&0&1&-1&-1&0 \\ 0&-2&0&1&2&-2&-1&0 \\ 0&-2&1&0&3&-3&-1&0 \\
        -1&-3&1&1&4&-4&-2&0 \\ 0&-2&0&1&3&-3&-2&0 \\
        0&-1&0&1&2&-3&-1&0 \\ 0&0&0&0&2&-2&-1&0 \\ 0&0&0&0&1&-1&-1&1 \end{array}\right)$,
  $\left(\begin{array}{cccccccc} -1&0&0&2&-2&0&0&0 \\ -1&0&1&2&-3&1&-1&0 \\ -2&0&1&3&-3&0&-1&0 \\
        -2&0&2&4&-5&0&-1&0 \\ -2&0&2&3&-4&0&-1&0 \\ -1&1&1&2&-3&0&-1&0 \\
        0&0&0&2&-2&0&-1&0 \\ 0&0&0&1&-1&0&0&-1 \end{array}\right)$,
  $\left(\begin{array}{cccccccc} 0&-1&1&0&1&-1&-1&0 \\ 0&-2&0&1&2&-2&-1&0 \\ 0&-2&1&0&3&-3&-1&0 \\
        -1&-3&1&1&4&-4&-2&0 \\ 0&-2&0&1&3&-3&-2&0 \\
        0&-1&0&1&2&-3&-1&0 \\ 0&0&0&0&2&-2&-1&0 \\ 0&0&0&0&1&-1&0&-1 \end{array}\right)$,
  $\left(\begin{array}{cccccccc} -1&-1&0&2&-2&1&0&0 \\ -1&-1&1&2&-3&2&-1&0 \\
        -2&-2&1&3&-3&2&-1&0 \\ -2&-3&2&4&-5&3&-1&0 \\
        -2&-3&2&3&-4&3&-1&0 \\ -1&-2&1&2&-3&3&-1&0 \\
        0&-2&0&2&-2&2&-1&0 \\ 0&-1&0&1&-1&1&-1&1 \end{array}\right)$,
  $\left(\begin{array}{cccccccc} 0&-2&1&0&1&0&-1&0 \\ 0&-3&0&1&2&-1&-1&0 \\ 0&-4&1&0&3&-1&-1&0 \\
        -1&-6&1&1&4&-1&-2&0 \\ 0&-5&0&1&3&0&-2&0 \\ 0&-4&0&1&2&0&-1&0 \\
        0&-2&0&0&2&0&-1&0 \\ 0&-1&0&0&1&0&-1&1 \end{array}\right)$,
  $\left(\begin{array}{cccccccc} -1&-1&0&2&-2&1&0&0 \\ -1&-1&1&2&-3&2&-1&0 \\ -2&-2&1&3&-3&2&-1&0 \\
        -2&-3&2&4&-5&3&-1&0 \\ -2&-3&2&3&-4&3&-1&0 \\
        -1&-2&1&2&-3&3&-1&0 \\ 0&-2&0&2&-2&2&-1&0 \\
        0&-1&0&1&-1&1&0&-1 \end{array}\right)$,
  $\left(\begin{array}{cccccccc} 0&-2&1&0&1&0&-1&0 \\ 0&-3&0&1&2&-1&-1&0 \\ 0&-4&1&0&3&-1&-1&0 \\
        -1&-6&1&1&4&-1&-2&0 \\ 0&-5&0&1&3&0&-2&0 \\ 0&-4&0&1&2&0&-1&0 \\
        0&-2&0&0&2&0&-1&0 \\ 0&-1&0&0&1&0&0&-1 \end{array}\right)$,
  $\left(\begin{array}{cccccccc} -1&0&0&2&-2&0&0&0 \\ 0&-1&0&3&-3&0&0&0 \\ -1&-1&0&4&-3&-1&0&0 \\
        -1&-1&1&5&-5&-1&0&0 \\ -1&-1&1&4&-4&-1&0&0 \\
        0&0&0&3&-3&-1&0&0 \\ 0&0&0&2&-2&0&-1&0 \\ 0&0&0&1&-1&0&-1&1 \end{array}\right)$,
  $\left(\begin{array}{cccccccc} 0&-1&1&0&1&-1&-1&0 \\ -1&-1&1&0&2&-1&-2&0 \\
        -1&-1&2&-1&3&-2&-2&0 \\ -2&-2&2&0&4&-3&-3&0 \\
        -1&-1&1&0&3&-2&-3&0 \\ -1&0&1&0&2&-2&-2&0 \\ 0&0&0&0&2&-2&-1&0 \\
        0&0&0&0&1&-1&-1&1 \end{array}\right)$,
  $\left(\begin{array}{cccccccc} -1&0&0&2&-2&0&0&0 \\ 0&-1&0&3&-3&0&0&0 \\ -1&-1&0&4&-3&-1&0&0 \\
        -1&-1&1&5&-5&-1&0&0 \\ -1&-1&1&4&-4&-1&0&0 \\
        0&0&0&3&-3&-1&0&0 \\ 0&0&0&2&-2&0&-1&0 \\ 0&0&0&1&-1&0&0&-1 \end{array}\right)$,
 $\left(\begin{array}{cccccccc} 0&-1&1&0&1&-1&-1&0 \\ -1&-1&1&0&2&-1&-2&0 \\
        -1&-1&2&-1&3&-2&-2&0 \\ -2&-2&2&0&4&-3&-3&0 \\
        -1&-1&1&0&3&-2&-3&0 \\ -1&0&1&0&2&-2&-2&0 \\ 0&0&0&0&2&-2&-1&0 \\
        0&0&0&0&1&-1&0&-1 \end{array}\right)$,
  $\left(\begin{array}{cccccccc} -1&-1&0&2&-2&1&0&0 \\ 0&-2&0&3&-3&1&0&0 \\ -1&-3&0&4&-3&1&0&0 \\
        -1&-4&1&5&-5&2&0&0 \\ -1&-4&1&4&-4&2&0&0 \\ 0&-3&0&3&-3&2&0&0 \\
        0&-2&0&2&-2&2&-1&0 \\ 0&-1&0&1&-1&1&-1&1 \end{array}\right)$,
  $\left(\begin{array}{cccccccc} 0&-2&1&0&1&0&-1&0 \\ -1&-2&1&0&2&0&-2&0 \\ -1&-3&2&-1&3&0&-2&0 \\
        -2&-5&2&0&4&0&-3&0 \\ -1&-4&1&0&3&1&-3&0 \\ -1&-3&1&0&2&1&-2&0 \\
        0&-2&0&0&2&0&-1&0 \\ 0&-1&0&0&1&0&-1&1 \end{array}\right)$,
  $\left(\begin{array}{cccccccc} 0&0&1&1&-2&0&0&0 \\ -1&0&1&2&-3&1&-1&0 \\ -1&0&2&2&-3&0&-1&0 \\
        -2&0&2&4&-5&0&-1&0 \\ -2&0&2&3&-4&0&-1&0 \\ -1&1&1&2&-3&0&-1&0 \\
        0&0&0&2&-2&0&-1&0 \\ 0&0&0&1&-1&0&-1&1 \end{array}\right)$,
  $\left(\begin{array}{cccccccc} -1&-1&0&1&1&-1&-1&0 \\ 0&-2&0&1&2&-2&-1&0 \\ -1&-2&0&1&3&-3&-1&0 \\
        -1&-3&1&1&4&-4&-2&0 \\ 0&-2&0&1&3&-3&-2&0 \\
        0&-1&0&1&2&-3&-1&0 \\ 0&0&0&0&2&-2&-1&0 \\ 0&0&0&0&1&-1&-1&1 \end{array}\right)$,
  $\left(\begin{array}{cccccccc} 0&0&1&1&-2&0&0&0 \\ -1&0&1&2&-3&1&-1&0 \\ -1&0&2&2&-3&0&-1&0 \\
        -2&0&2&4&-5&0&-1&0 \\ -2&0&2&3&-4&0&-1&0 \\ -1&1&1&2&-3&0&-1&0 \\
        0&0&0&2&-2&0&-1&0 \\ 0&0&0&1&-1&0&0&-1 \end{array}\right)$,
  $\left(\begin{array}{cccccccc} -1&-1&0&1&1&-1&-1&0 \\ 0&-2&0&1&2&-2&-1&0 \\ -1&-2&0&1&3&-3&-1&0 \\
        -1&-3&1&1&4&-4&-2&0 \\ 0&-2&0&1&3&-3&-2&0 \\
        0&-1&0&1&2&-3&-1&0 \\ 0&0&0&0&2&-2&-1&0 \\ 0&0&0&0&1&-1&0&-1 \end{array}\right)$,
  $\left(\begin{array}{cccccccc} 0&-1&1&1&-2&1&0&0 \\ -1&-1&1&2&-3&2&-1&0 \\
        -1&-2&2&2&-3&2&-1&0 \\ -2&-3&2&4&-5&3&-1&0 \\
        -2&-3&2&3&-4&3&-1&0 \\ -1&-2&1&2&-3&3&-1&0 \\
        0&-2&0&2&-2&2&-1&0 \\ 0&-1&0&1&-1&1&-1&1 \end{array}\right)$,
  $\left(\begin{array}{cccccccc} -1&-2&0&1&1&0&-1&0 \\ 0&-3&0&1&2&-1&-1&0 \\ -1&-4&0&1&3&-1&-1&0 \\
        -1&-6&1&1&4&-1&-2&0 \\ 0&-5&0&1&3&0&-2&0 \\ 0&-4&0&1&2&0&-1&0 \\
        0&-2&0&0&2&0&-1&0 \\ 0&-1&0&0&1&0&-1&1 \end{array}\right)$,
  $\left(\begin{array}{cccccccc} 1&0&-1&2&-2&0&0&0 \\ 0&0&-1&3&-3&1&-1&0 \\ 1&0&-2&4&-3&0&-1&0 \\
        1&1&-2&5&-5&1&-2&0 \\ 0&1&-1&4&-4&0&-1&0 \\ 0&1&-1&3&-3&0&-1&0 \\
        0&0&0&2&-2&0&-1&0 \\ 0&0&0&1&-1&0&-1&1 \end{array}\right)$,
  $\left(\begin{array}{cccccccc} -1&0&1&0&1&-1&-1&0 \\ 0&-1&1&0&2&-2&-1&0 \\ -1&0&2&-1&3&-3&-1&0 \\
        -1&-1&2&0&4&-5&-1&0 \\ -1&-1&2&0&3&-4&-1&0 \\
        0&0&1&0&2&-3&-1&0 \\ 0&0&0&0&2&-2&-1&0 \\ 0&0&0&0&1&-1&-1&1 \end{array}\right)$,
  $\left(\begin{array}{cccccccc} 1&0&-1&2&-2&0&0&0 \\ 0&0&-1&3&-3&1&-1&0 \\ 1&0&-2&4&-3&0&-1&0 \\
        1&1&-2&5&-5&1&-2&0 \\ 0&1&-1&4&-4&0&-1&0 \\ 0&1&-1&3&-3&0&-1&0 \\
        0&0&0&2&-2&0&-1&0 \\ 0&0&0&1&-1&0&0&-1 \end{array}\right)$,
  $\left(\begin{array}{cccccccc} -1&0&1&0&1&-1&-1&0 \\ 0&-1&1&0&2&-2&-1&0 \\ -1&0&2&-1&3&-3&-1&0 \\
        -1&-1&2&0&4&-5&-1&0 \\ -1&-1&2&0&3&-4&-1&0 \\
        0&0&1&0&2&-3&-1&0 \\ 0&0&0&0&2&-2&-1&0 \\ 0&0&0&0&1&-1&0&-1 \end{array}\right)$,
   $\left(\begin{array}{cccccccc} 1&-1&-1&2&-2&1&0&0 \\ 0&-1&-1&3&-3&2&-1&0 \\
        1&-2&-2&4&-3&2&-1&0 \\ 1&-2&-2&5&-5&4&-2&0 \\
        0&-2&-1&4&-4&3&-1&0 \\ 0&-2&-1&3&-3&3&-1&0 \\
        0&-2&0&2&-2&2&-1&0 \\ 0&-1&0&1&-1&1&-1&1 \end{array}\right)$,
  $\left(\begin{array}{cccccccc} -1&-1&1&0&1&0&-1&0 \\ 0&-2&1&0&2&-1&-1&0 \\ -1&-2&2&-1&3&-1&-1&0 \\
        -1&-4&2&0&4&-2&-1&0 \\ -1&-4&2&0&3&-1&-1&0 \\
        0&-3&1&0&2&0&-1&0 \\ 0&-2&0&0&2&0&-1&0 \\ 0&-1&0&0&1&0&-1&1 \end{array}\right)$,
  $\left(\begin{array}{cccccccc} 0&1&0&1&-2&1&-1&0 \\ 0&0&-1&3&-3&1&-1&0 \\ 0&1&-1&3&-3&1&-2&0 \\
        1&1&-2&5&-5&1&-2&0 \\ 0&1&-1&4&-4&0&-1&0 \\ 0&1&-1&3&-3&0&-1&0 \\
        0&0&0&2&-2&0&-1&0 \\ 0&0&0&1&-1&0&-1&1 \end{array}\right)$,
  $\left(\begin{array}{cccccccc} 0&-1&0&1&1&-2&0&0 \\ 0&-1&1&0&2&-2&-1&0 \\ 0&-1&1&0&3&-4&0&0 \\
        -1&-1&2&0&4&-5&-1&0 \\ -1&-1&2&0&3&-4&-1&0 \\
        0&0&1&0&2&-3&-1&0 \\ 0&0&0&0&2&-2&-1&0 \\ 0&0&0&0&1&-1&-1&1 \end{array}\right)$,
  $\left(\begin{array}{cccccccc} 0&1&0&1&-2&1&-1&0 \\ 0&0&-1&3&-3&1&-1&0 \\ 0&1&-1&3&-3&1&-2&0 \\
        1&1&-2&5&-5&1&-2&0 \\ 0&1&-1&4&-4&0&-1&0 \\ 0&1&-1&3&-3&0&-1&0 \\
        0&0&0&2&-2&0&-1&0 \\ 0&0&0&1&-1&0&0&-1 \end{array}\right)$,
  $\left(\begin{array}{cccccccc} 0&-1&0&1&1&-2&0&0 \\ 0&-1&1&0&2&-2&-1&0 \\ 0&-1&1&0&3&-4&0&0 \\
        -1&-1&2&0&4&-5&-1&0 \\ -1&-1&2&0&3&-4&-1&0 \\
        0&0&1&0&2&-3&-1&0 \\ 0&0&0&0&2&-2&-1&0 \\ 0&0&0&0&1&-1&0&-1 \end{array}\right)$,
  $\left(\begin{array}{cccccccc} 0&0&0&1&-2&2&-1&0 \\ 0&-1&-1&3&-3&2&-1&0 \\
        0&-1&-1&3&-3&3&-2&0 \\ 1&-2&-2&5&-5&4&-2&0 \\
        0&-2&-1&4&-4&3&-1&0 \\ 0&-2&-1&3&-3&3&-1&0 \\
        0&-2&0&2&-2&2&-1&0 \\ 0&-1&0&1&-1&1&-1&1 \end{array}\right)$,
  $\left(\begin{array}{cccccccc} 0&-2&0&1&1&-1&0&0 \\ 0&-2&1&0&2&-1&-1&0 \\ 0&-3&1&0&3&-2&0&0 \\
        -1&-4&2&0&4&-2&-1&0 \\ -1&-4&2&0&3&-1&-1&0 \\
        0&-3&1&0&2&0&-1&0 \\ 0&-2&0&0&2&0&-1&0 \\ 0&-1&0&0&1&0&-1&1 \end{array}\right)$,
  $\left(\begin{array}{cccccccc} -1&-1&0&1&1&-1&0&-1 \\ -1&-1&1&0&2&-1&-1&-1 \\
        -2&-1&1&0&3&-2&0&-2 \\ -2&-2&2&0&4&-3&0&-3 \\
        -2&-2&2&0&3&-2&0&-3 \\ -1&-1&1&0&2&-1&0&-3 \\
        0&-1&0&0&2&-1&0&-2 \\ 0&-1&0&0&1&0&0&-1 \end{array}\right)$,
  $\left(\begin{array}{cccccccc} 0&0&1&-2&2&-1&1&-1 \\ 0&0&0&-2&3&-2&2&-1 \\ 0&0&1&-3&3&-2&3&-2 \\
        -1&0&1&-4&5&-3&4&-3 \\ 0&0&0&-3&4&-2&3&-3 \\
        0&0&0&-2&3&-2&3&-3 \\ 0&1&0&-2&2&-1&2&-2 \\ 0&0&0&-1&1&0&1&-1 \end{array}\right)$,
  $\left(\begin{array}{cccccccc} -1&-1&0&1&1&-1&-1&1 \\ -1&-1&1&0&2&-1&-2&1 \\
        -2&-1&1&0&3&-2&-2&2 \\ -2&-2&2&0&4&-3&-3&3 \\
        -2&-2&2&0&3&-2&-3&3 \\ -1&-1&1&0&2&-1&-3&3 \\
        0&-1&0&0&2&-1&-2&2 \\ 0&-1&0&0&1&0&-1&1 \end{array}\right)$,
  $\left(\begin{array}{cccccccc} 0&0&1&-2&2&-1&0&1 \\ 0&0&0&-2&3&-2&1&1 \\ 0&0&1&-3&3&-2&1&2 \\
        -1&0&1&-4&5&-3&1&3 \\ 0&0&0&-3&4&-2&0&3 \\ 0&0&0&-2&3&-2&0&3 \\
        0&1&0&-2&2&-1&0&2 \\ 0&0&0&-1&1&0&0&1 \end{array}\right)$,
  $\left(\begin{array}{cccccccc} -1&-1&0&1&1&-1&0&-1 \\ 0&-2&0&1&2&-2&0&-1 \\ -1&-2&0&1&3&-3&1&-2 \\
        -1&-3&1&1&4&-4&1&-3 \\ -1&-3&1&1&3&-3&1&-3 \\
        0&-2&0&1&2&-2&1&-3 \\ 0&-1&0&0&2&-1&0&-2 \\ 0&-1&0&0&1&0&0&-1 \end{array}\right)$,
  $\left(\begin{array}{cccccccc} 0&0&1&-2&2&-1&1&-1 \\ -1&1&1&-3&3&-1&1&-1 \\
        -1&1&2&-4&3&-1&2&-2 \\ -2&1&2&-5&5&-2&3&-3 \\
        -1&1&1&-4&4&-1&2&-3 \\ -1&1&1&-3&3&-1&2&-3 \\
        0&1&0&-2&2&-1&2&-2 \\ 0&0&0&-1&1&0&1&-1 \end{array}\right)$,
  $\left(\begin{array}{cccccccc} -1&-1&0&1&1&-1&-1&1 \\ 0&-2&0&1&2&-2&-1&1 \\ -1&-2&0&1&3&-3&-1&2 \\
        -1&-3&1&1&4&-4&-2&3 \\ -1&-3&1&1&3&-3&-2&3 \\
        0&-2&0&1&2&-2&-2&3 \\ 0&-1&0&0&2&-1&-2&2 \\ 0&-1&0&0&1&0&-1&1 \end{array}\right)$,
  $\left(\begin{array}{cccccccc} 0&0&1&-2&2&-1&0&1 \\ -1&1&1&-3&3&-1&0&1 \\ -1&1&2&-4&3&-1&0&2 \\
        -2&1&2&-5&5&-2&0&3 \\ -1&1&1&-4&4&-1&-1&3 \\
        -1&1&1&-3&3&-1&-1&3 \\ 0&1&0&-2&2&-1&0&2 \\ 0&0&0&-1&1&0&0&1 \end{array}\right)$,
  $\left(\begin{array}{cccccccc} 0&-1&1&0&1&-1&0&-1 \\ -1&-1&1&0&2&-1&-1&-1 \\
        -1&-1&2&-1&3&-2&0&-2 \\ -2&-2&2&0&4&-3&0&-3 \\
        -2&-2&2&0&3&-2&0&-3 \\ -1&-1&1&0&2&-1&0&-3 \\
        0&-1&0&0&2&-1&0&-2 \\ 0&-1&0&0&1&0&0&-1 \end{array}\right)$,
  $\left(\begin{array}{cccccccc} -1&0&0&-1&2&-1&1&-1 \\ 0&0&0&-2&3&-2&2&-1 \\ -1&0&0&-2&3&-2&3&-2 \\
        -1&0&1&-4&5&-3&4&-3 \\ 0&0&0&-3&4&-2&3&-3 \\
        0&0&0&-2&3&-2&3&-3 \\ 0&1&0&-2&2&-1&2&-2 \\ 0&0&0&-1&1&0&1&-1 \end{array}\right)$,
 $\left(\begin{array}{cccccccc} 0&-1&1&0&1&-1&-1&1 \\ -1&-1&1&0&2&-1&-2&1 \\
        -1&-1&2&-1&3&-2&-2&2 \\ -2&-2&2&0&4&-3&-3&3 \\
        -2&-2&2&0&3&-2&-3&3 \\ -1&-1&1&0&2&-1&-3&3 \\
        0&-1&0&0&2&-1&-2&2 \\ 0&-1&0&0&1&0&-1&1 \end{array}\right)$,
  $\left(\begin{array}{cccccccc} -1&0&0&-1&2&-1&0&1 \\ 0&0&0&-2&3&-2&1&1 \\ -1&0&0&-2&3&-2&1&2 \\
        -1&0&1&-4&5&-3&1&3 \\ 0&0&0&-3&4&-2&0&3 \\ 0&0&0&-2&3&-2&0&3 \\
        0&1&0&-2&2&-1&0&2 \\ 0&0&0&-1&1&0&0&1 \end{array}\right)$,
  $\left(\begin{array}{cccccccc} 1&-1&-1&1&1&-1&0&-1 \\ 0&-1&-1&1&2&-1&-1&-1 \\
        1&-1&-2&1&3&-2&0&-2 \\ 1&-1&-2&1&4&-2&-1&-3 \\
        0&-1&-1&1&3&-2&0&-3 \\ 0&-1&-1&1&2&-1&0&-3 \\
        0&-1&0&0&2&-1&0&-2 \\ 0&-1&0&0&1&0&0&-1 \end{array}\right)$,
  $\left(\begin{array}{cccccccc} -1&1&1&-2&2&-1&1&-1 \\ 0&1&1&-3&3&-2&2&-1 \\ -1&2&2&-4&3&-2&3&-2 \\
        -1&2&2&-5&5&-4&5&-3 \\ -1&1&2&-4&4&-3&4&-3 \\
        0&1&1&-3&3&-2&3&-3 \\ 0&1&0&-2&2&-1&2&-2 \\ 0&0&0&-1&1&0&1&-1 \end{array}\right)$,
  $\left(\begin{array}{cccccccc} 1&-1&-1&1&1&-1&-1&1 \\ 0&-1&-1&1&2&-1&-2&1 \\
        1&-1&-2&1&3&-2&-2&2 \\ 1&-1&-2&1&4&-2&-4&3 \\
        0&-1&-1&1&3&-2&-3&3 \\ 0&-1&-1&1&2&-1&-3&3 \\
        0&-1&0&0&2&-1&-2&2 \\ 0&-1&0&0&1&0&-1&1 \end{array}\right)$,
  $\left(\begin{array}{cccccccc} -1&1&1&-2&2&-1&0&1 \\ 0&1&1&-3&3&-2&1&1 \\ -1&2&2&-4&3&-2&1&2 \\
        -1&2&2&-5&5&-4&2&3 \\ -1&1&2&-4&4&-3&1&3 \\ 0&1&1&-3&3&-2&0&3 \\
        0&1&0&-2&2&-1&0&2 \\ 0&0&0&-1&1&0&0&1 \end{array}\right)$,
  $\left(\begin{array}{cccccccc} 0&0&0&0&1&0&-1&-1 \\ 0&-1&-1&1&2&-1&-1&-1 \\ 0&0&-1&0&3&-1&-1&-2 \\
        1&-1&-2&1&4&-2&-1&-3 \\ 0&-1&-1&1&3&-2&0&-3 \\
        0&-1&-1&1&2&-1&0&-3 \\ 0&-1&0&0&2&-1&0&-2 \\ 0&-1&0&0&1&0&0&-1 \end{array}\right)$,
  $\left(\begin{array}{cccccccc} 0&0&0&-1&2&-2&2&-1 \\ 0&1&1&-3&3&-2&2&-1 \\
        0&1&1&-3&3&-3&4&-2 \\ -1&2&2&-5&5&-4&5&-3 \\
        -1&1&2&-4&4&-3&4&-3 \\ 0&1&1&-3&3&-2&3&-3 \\
        0&1&0&-2&2&-1&2&-2 \\ 0&0&0&-1&1&0&1&-1 \end{array}\right)$,
  $\left(\begin{array}{cccccccc} 0&0&0&0&1&0&-2&1 \\ 0&-1&-1&1&2&-1&-2&1 \\ 0&0&-1&0&3&-1&-3&2 \\
        1&-1&-2&1&4&-2&-4&3 \\ 0&-1&-1&1&3&-2&-3&3 \\
        0&-1&-1&1&2&-1&-3&3 \\ 0&-1&0&0&2&-1&-2&2 \\ 0&-1&0&0&1&0&-1&1 \end{array}\right)$,
  $\left(\begin{array}{cccccccc} 0&0&0&-1&2&-2&1&1 \\ 0&1&1&-3&3&-2&1&1 \\ 0&1&1&-3&3&-3&2&2 \\
        -1&2&2&-5&5&-4&2&3 \\ -1&1&2&-4&4&-3&1&3 \\ 0&1&1&-3&3&-2&0&3 \\
        0&1&0&-2&2&-1&0&2 \\ 0&0&0&-1&1&0&0&1 \end{array}\right)$.

The character table of $G^{s_{47}}$:\\

\noindent \begin{tabular}{c|cccccccccccccccccccccccccccccc}
  & & & & & & & & & & 10& & & & & & & & & &20 & & & & & & & & & &30\\\hline
$\chi_{47}^ {(1)}$&1&1&1&1&1&1&1&1&1&1&1&1&1&1&1&1&1&1&1&1
&1&1&1&1&1&1&1&1&1&1

\\    $\chi_{47}^ {(2)}$&1&1&1&1&-1&-1&-1&-1&1&1&1&1&-1&-1&-1&-1&1&1&1&1 &-1&-1&-1
&-1&1&1&1&1&-1&-1

\\    $\chi_{47}^ {(3)}$&1&1&1&1&-1&-1&-1&-1&1&1&1&1&-1&-1&-1&-1&1&1&1&1 &-1&-1&-1
&-1&1&1&1&1&-1&-1

\\    $\chi_{47}^ {(4)}$&1&1&1&1&-1&-1&-1&-1&1&1&1&1&-1&-1&-1&-1&1&1&1&1 &-1&-1&-1
&-1&1&1&1&1&-1&-1

\\    $\chi_{47}^ {(5)}$&1&1&1&1&-1&-1&-1&-1&1&1&1&1&-1&-1&-1&-1&1&1&1&1 &-1&-1&-1
&-1&1&1&1&1&-1&-1

\\    $\chi_{47}^ {(6)}$&1&1&1&1&-1&-1&-1&-1&1&1&1&1&-1&-1&-1&-1&1&1&1&1 &-1&-1&-1
&-1&1&1&1&1&-1&-1

\\    $\chi_{47}^ {(7)}$&1&1&1&1&-1&-1&-1&-1&1&1&1&1&-1&-1&-1&-1&1&1&1&-1&-1&-1
&-1&1&1&1&1&-1&-1

\\
\end {tabular}

\noindent \begin{tabular}{c|cccccccccccccccccccccccccccccc}
  & & & & & & & & & & 10& & & & & & & & & &20 & & & & & & & & & &30\\\hline

\\    $\chi_{47}^ {({46})}$&1&1&-1&-1&-1&-1&1&1&1&1&-1&-1&1&1&-1&-1&-1&-1&1&1&1&1&-1
&-1&1&1&-1&-1&-1&-1

\\    $\chi_{47}^ {({47})}$&1&1&-1&-1&-1&-1&1&1&1&1&-1&-1&1&1&-1&-1&-1&-1&1&1&1&1&-1
&-1&1&1&-1&-1&-1&-1

\\    $\chi_{47}^ {({48})}$&1&1&-1&-1&-1&-1&1&1&1&1&-1&-1&1&1&-1&-1&-1&-1&1&1&1&1&-1
&-1&1&1&-1&-1&-1&-1

\\    $\chi_{47}^ {({49})}$&1&1&-1&-1&1&1&-1&-1&1&1&-1&-1&-1&-1&1&1&-1&-1&1&1&-1
&-1&1&1&1&1&-1&-1&1&1

\\    $\chi_{47}^ {({50})}$&1&1&-1&-1&1&1&-1&-1&1&1&-1&-1&-1&-1&1&1&-1&-1&1&1&-1
&-1&1&1&1&1&-1&-1&1&1

\\    $\chi_{47}^ {({51})}$&1&1&-1&-1&1&1&-1&-1&1&1&-1&-1&-1&-1&1&1&-1&-1&1&1&-1
&-1&1&1&1&1&-1&-1&1&1

\\
\end {tabular}

${}_{
 \!\!\!\!\!\!\!\!\!\!\!\!\!\!\!\!\!\!_{
\!\!\!\!\!\!\!_{ \!\!\!\!\!\!\!  _{  \small
 \noindent
\begin{tabular}{c|cccccccccccccccccccccccccccccc}
  & & & & & & & & & & 10& & & & & & & & & &20 & & & & & & & & & &30\\\hline
\\    $\chi_{47}^ {({91})}$&2&2&2&2&.&.&.&.&-2&-2&-2&-2&-2&-2&-2&-2&.&.&.&.
&2&2&2&2&2&2&2&2&.&.

\\    $\chi_{47}^ {({92})}$&2&2&2&2&.&.&.&.&-2&-2&-2&-2&-2&-2&-2&-2&.&.&.&.
&2&2&2&2&2&2&2&2&.&.

\\    $\chi_{47}^ {({93})}$&2&2&2&2&.&.&.&.&-2&-2&-2&-2&2&2&2&2&.&.&.&.&-2&-2&-2
&-2&2&2&2&2&.&.

\\    $\chi_{47}^ {({94})}$&2&2&2&2&.&.&.&.&-2&-2&-2&-2&2&2&2&2&.&.&.&.&-2&-2&-2&-2&2&2&2
&2&.&.

\\    $\chi_{47}^ {({95})}$&2&2&2&2&.&.&.&.&-2&-2&-2&-2&2&2&2&2&.&.&.&.&-2&-2&-2
&-2&2&2&2&2&.&.

\\    $\chi_{47}^ {({96})}$&2&2&2&2&.&.&.&.&-2&-2&-2&-2&2&2&2&2&.&.&.&.&-2&-2&-2
&-2&2&2&2&2&.&.

\\
\end {tabular}
 }}}}$

${}_{
 \!\!\!\!\!\!\!_{
\!\!\!\!\!\!\!_{ \!\!\!\!\!\!\!  _{  \small
 \noindent
\begin{tabular}{c|cccccccccccccccccccccccccccccc}
  & & & & & & & & & & 10& & & & & & & & & &20 & & & & & & & & & &30\\\hline
\\    $\chi_{47}^ {({136})}$&2&2&2&2&-2&-2&-2&-2&2&2&2&2&.&.&.&.&.&.&.&.&.&.&.&.&-2&-2&-2
&-2&2&2

\\    $\chi_{47}^ {({137})}$&2&2&-2&-2&.&.&.&.&-2&-2&2&2&-2&-2&2&2&.&.&.&.&2&2&-2&-2&2&2
&-2&-2&.&.

\\    $\chi_{47}^ {({138})}$&2&2&-2&-2&.&.&.&.&-2&-2&2&2&-2&-2&2&2&.&.&.&.&2&2&-2&-2&2&2
&-2&-2&.&.

\\    $\chi_{47}^ {({139})}$&2&2&-2&-2&.&.&.&.&-2&-2&2&2&-2&-2&2&2&.&.&.&.&2&2&-2&-2&2&2
&-2&-2&.&.

\\    $\chi_{47}^ {({140})}$&2&2&-2&-2&.&.&.&.&-2&-2&2&2&-2&-2&2&2&.&.&.&.&2&2&-2&-2&2&2
&-2&-2&.&.

\\    $\chi_{47}^ {({141})}$&2&2&-2&-2&.&.&.&.&-2&-2&2&2&2&2&-2&-2&.&.&.&.&-2&-2&2&2&2&2
&-2&-2&.&.

\\
\end{tabular}
 }}}}$

\noindent \begin{tabular}{c|cccccccccccccccccccccccccccc}
  & & & & & & & & & &40 & & & & & & & & & & 50& & & & & & & &\\\hline
$\chi_{47}^ {(1)}$&1&1&1&1&1&1&1&1&1&1&1&1&1&1&1&1&1&1&1&
1&1&1&1&1&1&1&1&1

\\    $\chi_{47}^ {(2)}$&-1&-1&1&1&1&1&1&1&-1&-1&1&1&-1&-1&1&1&-1&-1&1&1&-1
&-1&1&1&-1&-1&-1&-1

\\    $\chi_{47}^ {(3)}$&-1&-1&1&1&1&1&1&1&-1&-1&1&1&-1&-1&1&1&-1&-1&1&1&-1
&-1&1&1&-1&-1&-1&-1

\\    $\chi_{47}^ {(4)}$&-1&-1&1&1&1&1&1&1&-1&-1&1&1&-1&-1&1&1&-1&-1&1&1&-1
&-1&1&1&-1&-1&-1&-1

\\    $\chi_{47}^ {(5)}$&-1&-1&1&1&1&1&1&1&-1&-1&1&1&-1&-1&1&1&-1&-1&1&1&-1
&-1&1&1&-1&-1&-1&-1

\\    $\chi_{47}^ {(6)}$&-1&-1&1&1&1&1&1&1&-1&-1&1&1&-1&-1&1&1&-1&-1&1&1&-1
&-1&1&1&1&1&1&1

\\    $\chi_{47}^ {(7)}$&-1&-1&1&1&1&1&1&1&-1&-1&1&1&-1&-1&1&1&-1&-1&1&1&-1
&-1&1&1&1&1&1&1

\\
\end{tabular}

\noindent \begin{tabular}{c|cccccccccccccccccccccccccccc}
  & & & & & & & & & &40 & & & & & & & & & & 50& & & & & & & &\\\hline
\\    $\chi_{47}^ {({46})}$&1&1&1&1&-1&-1&-1&1&1&-1&-1&1&-1&1&1&-1&-1&1&-1&1&1&-1
&-1&1&-D&-D&D&D

\\    $\chi_{47}^ {({47})}$&1&1&1&1&-1&-1&-1&1&1&-1&-1&1&-1&1&1&-1&-1&1&-1&1&1&-1
&-1&1&D&D&-D&-D

\\    $\chi_{47}^ {({48})}$&1&1&1&1&-1&-1&-1&1&1&-1&-1&1&-1&1&1&-1&-1&1&-1&1&1&-1
&-1&1&-D&-D&D&D

\\    $\chi_{47}^ {({49})}$&-1&-1&1&1&-1&-1&-1&1&-1&1&-1&1&1&-1&1&-1&1&-1&-1&1&-1&1
&-1&1&D&D&-D&-D

\\    $\chi_{47}^ {({50})}$&-1&-1&1&1&-1&-1&-1&1&-1&1&-1&1&1&-1&1&-1&1&-1&-1&1&-1&1
&-1&1&-D&-D&D&D

\\    $\chi_{47}^ {({51})}$&-1&-1&1&1&-1&-1&-1&1&-1&1&-1&1&1&-1&1&-1&1&-1&-1&1&-1&1
&-1&1&D&D&-D&-D

\\
\end{tabular}

\noindent \begin{tabular}{c|cccccccccccccccccccccccccccc}
  & & & & & & & & & &40 & & & & & & & & & & 50& & & & & & & &\\\hline
\\    $\chi_{47}^ {({91})}$&.&.&-2&-2&-2&-2&2&2&.&.&-2&-2&-2&-2&.&.&2&2&2&2&.&.
&-2&-2&2&2&2&2

\\    $\chi_{47}^ {({92})}$&.&.&-2&-2&-2&-2&2&2&.&.&-2&-2&-2&-2&.&.&2&2&2&2&.&.
&-2&-2&-2&-2&-2&-2

\\    $\chi_{47}^ {({93})}$&.&.&-2&-2&-2&-2&2&2&.&.&-2&-2&2&2&.&.&-2&-2&2&2&.&.
&-2&-2&2&2&2&2

\\    $\chi_{47}^ {({94})}$&.&.&-2&-2&-2&-2&2&2&.&.&-2&-2&2&2&.&.&-2&-2&2&2&.&.
&-2&-2&-2&-2&-2&-2

\\    $\chi_{47}^ {({95})}$&.&.&-2&-2&-2&-2&2&2&.&.&-2&-2&2&2&.&.&-2&-2&2&2&.&.
&-2&-2&2&2&2&2

\\    $\chi_{47}^ {({96})}$&.&.&-2&-2&-2&-2&2&2&.&.&-2&-2&2&2&.&.&-2&-2&2&2&.&.
&-2&-2&-2&-2&-2&-2

\\
\end{tabular}%

\noindent \noindent where A = -2*E(4)
  = -2*ER(-1) = -2i,
B = -4*E(4)
  = -4*ER(-1) = -4i,
C = -8*E(4)
  = -8*ER(-1) = -8i,
D = -E(4)
  = -ER(-1) = -i,
E = -1-E(4)
  = -1-ER(-1) = -1-i,
F = -2-2*E(4)
  = -2-2*ER(-1) = -2-2i.

The generators of $G^{s_{48}}$ are:\\
$\left(%
\right)$.

The representatives of conjugacy classes of   $G^{s_{48}}$ are:\\
$\left(%
\right)$.

The character table of $G^{s_{48}}$:\\
\noindent %

A = E(3)$^2$
  = (-1-ER(-3))/2 = -1-b3,
B = 2*E(3)$^2$
  = -1-ER(-3) = -1-i3,
C = 3*E(3)$^2$
  = (-3-3*ER(-3))/2 = -3-3b3,
D = 4*E(3)$^2$
  = -2-2*ER(-3) = -2-2i3,
E = 6*E(3)$^2$
  = -3-3*ER(-3) = -3-3i3.

The generators of $G^{s_{49}}$ are:\\

$\left(%
\right)$.

The representatives of conjugacy classes of   $G^{s_{49}}$ are:\\

$\left(%
\right)$.

The character table of $G^{s_{49}}$:\\
${}_{
 \!\!\!\!\!\!\!_{
\!\!\!\!\!\!\!_{ \!\!\!\!\!\!\!  _{  \small
 \noindent
%

\noindent \noindent where A = E(4)
  = ER(-1) = i,
B = E(8)$^3$,C = 2*E(4)
  = 2*ER(-1) = 2i,
D = 2*E(8).

The generators of $G^{s_{50}}$ are:\\
$\left(%
\right)$.

The character table of $G^{s_{50}}$:\\
${}_{
 \!\!\!\!\!\!\!_{
\!\!\!\!\!\!\!_{ \!\!\!\!\!\!\!  _{  \small
 \noindent
%

\noindent \noindent where A = -2*E(4)
  = -2*ER(-1) = -2i,
B = -4*E(4)
  = -4*ER(-1) = -4i,
C = -6*E(4)
  = -6*ER(-1) = -6i,
D = E(4)
  = ER(-1) = i,
E = -1-E(4)
  = -1-ER(-1) = -1-i.

The generators of $G^{s_{51}}$ are:\\
$\left(%
\right)$.

The character table of $G^{s_{51}}$:\\
\noindent %

\noindent \noindent where A = -E(4)
  = -ER(-1) = -i,
B = 2*E(4)
  = 2*ER(-1) = 2i,
C = -3*E(4)
  = -3*ER(-1) = -3i.

The generators of $G^{s_{52}}$ are:\\
$\left(%
\right)$.

The character table of $G^{s_{52}}$:\\
${}_{
 \!\!\!\!\!\!\!_{
\!\!\!\!\!\!\!_{ \!\!\!\!\!\!\!  _{  \small
 \noindent
%

\noindent \noindent where A = E(3)$^2$
  = (-1-ER(-3))/2 = -1-b3,
B = 2*E(3)$^2$
  = -1-ER(-3) = -1-i3.

The generators of $G^{s_{53}}$ are:\\
$\left(%
\right)$.

The character table of $G^{s_{53}}$:\\
\noindent %

The generators of $G^{s_{54}}$ are:\\
$\left(%
\right)$.

The character table of $G^{s_{54}} =G$
:\\
${}_{
 \!\!\!\!\!\!\!_{
\!\!\!\!\!\!\!_{ \!\!\!\!\!\!\!  _{  \small
 \noindent
%

\noindent \noindent where A = -E(4)
  = -ER(-1) = -i,
B = E(3)
  = (-1+ER(-3))/2 = b3,
C = -E(12)$^7$.

The generators of $G^{s_{55}}$ are:\\
$\left(%
\right)$.

The character table of $G^{s_{55}}$:\\
\noindent %

\noindent \noindent where A = E(3)$^2$
  = (-1-ER(-3))/2 = -1-b3,
B = 2*E(3)$^2$
  = -1-ER(-3) = -1-i3,
C = 3*E(3)$^2$
  = (-3-3*ER(-3))/2 = -3-3b3,
D = 4*E(3)$^2$
  = -2-2*ER(-3) = -2-2i3,
E = 6*E(3)$^2$
  = -3-3*ER(-3) = -3-3i3.

The generators of $G^{s_{56}}$ are:\\
$\left(%
\right)$.

The character table of $G^{s_{56}}$:\\
\noindent %

\noindent \noindent where A = 3*E(3)$^2$
  = (-3-3*ER(-3))/2 = -3-3b3,
B = 6*E(3)$^2$
  = -3-3*ER(-3) = -3-3i3,
C = 9*E(3)$^2$
  = (-9-9*ER(-3))/2 = -9-9b3,
D = 12*E(3)$^2$
  = -6-6*ER(-3) = -6-6i3,
E = 18*E(3)$^2$
  = -9-9*ER(-3) = -9-9i3,
F = -E(3)$^2$
  = (1+ER(-3))/2 = 1+b3,
G = 2*E(3)$^2$
  = -1-ER(-3) = -1-i3,
H = 4*E(3)$^2$
  = -2-2*ER(-3) = -2-2i3,
I = -2*E(3)-E(3)$^2$
  = (3-ER(-3))/2 = 1-b3,
J = E(3)-E(3)$^2$
  = ER(-3) = i3,
K = -4*E(3)-2*E(3)$^2$
  = 3-ER(-3) = 3-i3,
L = 2*E(3)-2*E(3)$^2$
  = 2*ER(-3) = 2i3.

The generators of $G^{s_{57}}$ are:\\
$\left(%
\right)$.

The character table of $G^{s_{57}}$:\\
\noindent %

\noindent \noindent where A = 3*E(3)$^2$
  = (-3-3*ER(-3))/2 = -3-3b3,
B = 6*E(3)$^2$
  = -3-3*ER(-3) = -3-3i3,
C = 9*E(3)$^2$
  = (-9-9*ER(-3))/2 = -9-9b3,
D = 12*E(3)$^2$
  = -6-6*ER(-3) = -6-6i3,
E = 18*E(3)$^2$
  = -9-9*ER(-3) = -9-9i3,
F = -E(3)$^2$
  = (1+ER(-3))/2 = 1+b3,
G = 2*E(3)$^2$
  = -1-ER(-3) = -1-i3,
H = 4*E(3)$^2$
  = -2-2*ER(-3) = -2-2i3,
I = -E(3)+E(3)$^2$
  = -ER(-3) = -i3,
J = 2*E(3)+E(3)$^2$
  = (-3+ER(-3))/2 = -1+b3,
K = -2*E(3)+2*E(3)$^2$
  = -2*ER(-3) = -2i3,
L = 4*E(3)+2*E(3)$^2$
  = -3+ER(-3) = -3+i3.

The generators of $G^{s_{58}}$ are:\\
$\left(%
\right)$.

The character table of $G^{s_{58}}$:\\
\noindent %

\noindent \noindent where A = -E(3)$^2$
  = (1+ER(-3))/2 = 1+b3,
B = 2*E(3)$^2$
  = -1-ER(-3) = -1-i3,
C = 4*E(3)$^2$
  = -2-2*ER(-3) = -2-2i3.

The generators of $G^{s_{59}}$ are:\\
$\left(%
\right)$.

The character table of $G^{s_{59}}$:\\
${}_{
 \!\!\!\!\!\!\!_{
\!\!\!\!\!\!\!_{ \!\!\!\!\!\!\!  _{  \small
 \noindent
%

\noindent \noindent where A = E(4)
  = ER(-1) = i,
B = E(8)$^3$.

The generators of $G^{s_{60}}$ are:\\
$\left(%
\right)$.

The character table of $G^{s_{60}}$:\\
${}_{
 \!\!\!\!\!\!\!_{
\!\!\!\!\!\!\!_{ \!\!\!\!\!\!\!  _{  \small
 \noindent
%

\noindent \noindent where A = E(3)$^2$
  = (-1-ER(-3))/2 = -1-b3.

The generators of $G^{s_{61}}$ are:\\
$\left(%
\right)$.

The character table of $G^{s_{61}}$:\\
${}_{
 \!\!\!\!\!\!\!_{
\!\!\!\!\!\!\!_{ \!\!\!\!\!\!\!  _{  \small
 \noindent
%

\noindent \noindent where A = E(3)
  = (-1+ER(-3))/2 = b3,
B = -2*E(3)
  = 1-ER(-3) = 1-i3.

The generators of $G^{s_{62}}$ are:\\
$\left(%
\right)$.

The character table of $G^{s_{62}}$:\\
\noindent %

The generators of $G^{s_{63}}$ are:\\
$\left(%
\right)$.

The character table of $G^{s_{63}}$:\\
\noindent %

\noindent \noindent where A = E(3)$^2$
  = (-1-ER(-3))/2 = -1-b3,
B = 2*E(3)$^2$
  = -1-ER(-3) = -1-i3,
C = 4*E(3)$^2$
  = -2-2*ER(-3) = -2-2i3,
D = 8*E(3)$^2$
  = -4-4*ER(-3) = -4-4i3,
E = 3*E(3)$^2$
  = (-3-3*ER(-3))/2 = -3-3b3,
F = E(3)-E(3)$^2$
  = ER(-3) = i3,
G = 2*E(3)-2*E(3)$^2$
  = 2*ER(-3) = 2i3,
H = -2*E(3)-E(3)$^2$
  = (3-ER(-3))/2 = 1-b3,
I = -4*E(3)-2*E(3)$^2$
  = 3-ER(-3) = 3-i3.

The generators of $G^{s_{64}}$ are:\\
$\left(%
\right)$.

The character table of $G^{s_{64}}$:\\
${}_{
 \!\!\!\!\!\!\!_{
\!\!\!\!\!\!\!_{ \!\!\!\!\!\!\!  _{  \small
 \noindent
%

\noindent \noindent where A = E(3)$^2$
  = (-1-ER(-3))/2 = -1-b3,
B = -2*E(3)$^2$
  = 1+ER(-3) = 1+i3,
C = 3*E(3)$^2$
  = (-3-3*ER(-3))/2 = -3-3b3.
\vskip 0.3cm

\noindent{\large\bf Acknowledgement}:\\ We would like to thank Prof.
N. Andruskiewitsch and Dr. F. Fantino for suggestions and help.


\vskip 0.3cm

\noindent{\large\bf Acknowledgement}: We would like to thank Prof.
N. Andruskiewitsch and Dr. F. Fantino for suggestions and help.

}
\end {document}